\documentclass[10pt,english]{article}
\usepackage{geometry}
\usepackage{fancyhdr}
\usepackage{verbatim}
\usepackage{color}
\usepackage{mathrsfs}
\usepackage{amsmath}
\usepackage{amsthm}
\usepackage{amssymb}
\usepackage{hyperref}
\makeatletter
\theoremstyle{plain}
\newtheorem{theorem}{\protect\theoremname}[section] 
\theoremstyle{definition}
\newtheorem{definition}[theorem]{\protect\definitionname}
\theoremstyle{plain}
\newtheorem{lemma}[theorem]{\protect\lemmaname}
\theoremstyle{remark}
\newtheorem{remark}[theorem]{\protect\remarkname}
\theoremstyle{plain}
\newtheorem{corollary}[theorem]{\protect\corollaryname}
\theoremstyle{plain}
\newtheorem{proposition}[theorem]{\protect\propositionname}
\theoremstyle{plain}

\theoremstyle{plain}

\@ifundefined{date}{}{\date{}}
\usepackage{fancyhdr}

\usepackage{babel}
\providecommand{\definitionname}{Definition}
\providecommand{\lemmaname}{Lemma}
\providecommand{\theoremname}{Theorem}
\providecommand{\corollaryname}{Corollary}
\providecommand{\remarkname}{Remark}
\providecommand{\propositionname}{Proposition}
\providecommand{\examplename}{Example}
\providecommand{\probelmname}{Problem}

\makeatother
\def\rn{{\mathbb R^n}}
\def\d{{\mathrm d}}
\def\D{{\mathcal D}}
\begin{document}
	\title
	{\bf\Large
	Precompactness in matrix weighted Bourgain-Morrey spaces
		\footnotetext{Jingshi Xu is supported by the National Natural Science Foundation of China (Grant No. 12161022) and the Science and Technology Project of Guangxi (Guike AD23023002).
		}
		}
	
	\date{}
	
	\author{Tengfei Bai\textsuperscript{a}, Jingshi Xu\textsuperscript{b,c,d}\footnote{Corresponding author, E-mail: jingshixu@126.com}  \\
		{\scriptsize \textsuperscript{a} College of Mathematics and Statistics, Hainan Normal University, Haikou, Hainan 571158,
			China}\\
		{\scriptsize  \textsuperscript{b} School of Mathematics and Computing Science, Guilin University of Electronic Technology, Guilin 541004, China} \\
		{\scriptsize  \textsuperscript{c} Center for Applied Mathematics of Guangxi (GUET), Guilin 541004, China}\\
		{\scriptsize  \textsuperscript{d} Guangxi Colleges and Universities Key Laboratory of Data Analysis and Computation, Guilin 541004, China}
	}
	
	\pagestyle{myheadings}\markboth{\footnotesize\rm\sc Tengfei Bai, Pengfei Guo,  Jingshi Xu}
	{\footnotesize\rm\sc  	The preduals of  Banach space valued Bourgain-Morrey  spaces}
	
\maketitle
\begin{abstract}
	In this paper, we introduce matrix weighted Bourgain-Morrey spaces and obtain two sufficient conditions for precompact sets in matrix weighted Bourgain-Morrey spaces. We prove that the dyadic average operator is bounded on some matrix weighted Bourgain-Morrey spaces. With this result, we obtain the necessity for precompact sets in some matrix weighted Bourgain-Morrey spaces. The results are new even for the unweighted Bourgain-Morrey spaces.
\end{abstract}
\textbf{Keywords} Kolmogorov-Riesz theorem, matrix weight, precompactness,  Bourgain-Morrey space.

\noindent \textbf{Mathematics Subject Classification} 
Primary 46B50, 46E40, 42B35, 46E30.

\section{Introduction}

In 1931, Kolmogorov  \cite{K31} first discovered the characterization of precompact sets in $L^p ([0,1])$ for $p \in (1,\infty)$.
After that, there are many criteria for compactness of sets in Lebesgue spaces, which are called  the Kolmogorov-Riesz compactness theorems on the Lebesgue spaces. More details and the history, we refer the reader to \cite{HH10}.

Inspired by \cite{HH10}, Clop-Cruz \cite{CC13}  obtained a compactness criterion in scalar weighted Lebesgue spaces $L^p (\omega)$ for $ 1<p <\infty  $ with a scalar wight $\omega$ in Muckenhoupt class $A_p $. In \cite{GZ20}, Guo and Zhao improved the result in \cite{CC13} and obtained a compactness criterion in $L^p(\omega)$ for $p \in  (0,\infty)$ with
locally integrable weight $\omega$. In \cite{LYZ23}, Liu, Yang and  Zhuo proved the Kolmogorov-Riesz compactness theorem  in matrix  weighted Lebesgue spaces $L^p (W)$ with $1< p <\infty $.

The theory of matrix weighted function spaces goes back to  \cite{WM58}.
Indeed, in  1958, Wiener and Masani \cite[Section 4]{WM58} studied the matrix weighted $L^2 (W) $ for the  prediction theory for multivariate stochastic processes.
In \cite{TV97}, Treil and Volberg introduced matrix class $\mathcal A_2$. Nazarov and Treil \cite{NT97} and Volberg \cite{V97} extended $\mathcal A_2$ to $\mathcal A_p$ with $ p \in (1,\infty)$.
In \cite{G03}, Goldberg showed that the matrix $\mathcal A _p$ condition
leads to $L^p$ boundedness of a Hardy-Littlewood maximal function and obtained the boundedness  of matrix weighted singular integral operators in  Lebesgue spaces $L^p , 1<p<\infty$.
In \cite{R03,R04,R02}, Roudenko introduced the matrix-weighted homogeneous Besov spaces $ \dot B ^{s,q} _p (W)$ and matrix-weighted sequence Besov spaces $\dot b^{s,q}_p (W) $ and showed their equivalence via $\varphi$-transform and wavelets.
In \cite{FR04}, Frazier and Roudenko introduced the matrix class  $\mathcal A _p, (0<p \le 1)$
and studied the continuous and discrete matrix-weighted Besov spaces $\dot B^{s,q}_{p} (W)$ and $\dot b^{s,q}_p (W) $ with $0< p \le 1$.
In \cite{FR21}, Frazier and Roudenko introduced the homogeneous matrix-weighted  Triebel-Lizorkin spaces $\dot F^{s,q} _p (W) $ for $s\in \mathbb R, 0<p<\infty , 0<q\le \infty $ and
obtained the Littlewood-Paley characterizations of matrix-weighted Lebesgue spaces $L^p (W)$ and matrix-weighted Sobolev space $L^p_k (W)$ for $k\in \mathbb N, 1<p<\infty$.
In \cite{V97},  Volberg also introduced an analogue  condition for the matrix class  $\mathcal A_\infty$.
In  \cite{BHYY23,BHYY23II,BHYY23III}, Bu, Hyt\"onen,  Yang,  Yuan  studied the matrix weighted Besov-type and Triebel-Lizorkin-type spaces. Specifically, they introduced a new concept of the $\mathcal A_p$-dimension $\tilde d$, which  is useful to the proof of main results of this paper.
In \cite{BHYYinfty}, Bu, Hyt\"onen,  Yang,  Yuan obtained several new characterizations of $\mathcal A_{p,\infty}$-matrix weights.
	In \cite{BHYY25}, Bu et al. studied the  inhomogeneous Besov-type and Triebel--Lizorkin-type spaces with the result in \cite{BHYYinfty}.
	In  \cite{BCYY25}, Bu et al. introduced  the matrix weighted Hardy  spaces and obtained characterizations of these spaces via maximal function, atom. As applications, they established the finite atomic characterization of matrix weighted Hardy  spaces and obtained a criterion on the boundedness of sublinear operators from matrix weighted Hardy  spaces to any $\gamma$-quasi-Banach space ($\gamma\in (0,1]$). The boundedness of Calder\'on-Zygmund operators on matrix weighted Hardy  spaces was also obtained.
In \cite{LYY24}, Li, Yang and Yuan introduced the
matrix-weighted Besov-Triebel-Lizorkin spaces with logarithmic smoothness and
characterize these spaces via Peetre-type maximal functions. In \cite{ZYZ24,ZSTYY23}, Zhao et
al. introduced (generalized grand) Besov-Bourgain-Morrey spaces and explored
various real-variable properties of these spaces, which are a bridge connecting Bourgain-Morrey spaces with amalgam-type spaces. Moreover, some real-variable properties and boundedness of classical operators were studied in their article.
For many other results on the matrix class  $\mathcal A_p$   and matrix weighted function spaces, we refer the reader to \cite{BYY23,CMR16,FR04,FR21,G03,R03,R04,WYZ22}.

Bourgain \cite{B91} introduced a special case of Bourgain-Morrey spaces to study the Stein-Tomas (Strichartz) estimate.
In \cite{HNSH23}, Hatano, Nogayama,  Sawano, and  Hakim researched the Bourgain-Morrey spaces from the viewpoints of harmonic analysis and functional analysis.
In \cite{HLY23},  Hu, Li and Yang introduced the  Triebel-Lizorkin-Bourgain-Morrey spaces which connect Bourgain-Morrey spaces and global Morrey spaces.

Motivated by above literature, we will  introduce matrix weighted Bourgain-Morrey spaces and research precompact sets in these spaces.
The paper is   organized as follows. In Section \ref{sec 2}, dyadic cubes, the matrix  class  $\mathcal A _p$, $\mathcal A _p$-dimension $\tilde d$, matrix weighted Bourgain-Morrey spaces are given. The first result lies in Section \ref{KRct}. Specifically, a sufficient  condition for totally bounded set   in  matrix weighted spaces $ M _p^{t,r} (W)$ with $ 1\le p < r <\infty $ or $1 \le p  \le t <r =\infty $ is obtained in Theorem \ref{K-R 1}. As a application, we get a criterion for totally bounded set in degenerate Bourgain-Morrey spaces with matrix weight. The second result (Theorem \ref{K-R 2}) is replacing the translation operator by the average operator. Note that
the  translation operator is not bounded on $L^p (W)$ and $M_p^{t,r} (W) $ in general.  We prove that dyadic average operator is  bounded on matrix weighted Bourgain-Morrey spaces with some conditions in Theorem \ref{tm ave BM}. Using this result, we obtain the Kolmogorov-Riesz compactness theorem in matrix weighted Bourgain-Morrey spaces. These results are new even for the unweighted Bourgain-Morrey spaces.

Throughout this paper, we let $c, C$ denote constants that are independent of the main parameters involved but whose value may differ from line to line.
Let $\mathbb N = \{1,2,3,\ldots \}$ and $\mathbb N_0 = \{0,1,2,3,\ldots \}$. Let $\mathbb Z$ be the set of integers.
Let $\chi_{E}$ be the characteristic function of the set $E\subset \rn $.
By $A\lesssim B$ we mean that $A\leq CB$ with some positive constant $C$ independent of appropriate quantities. By $ A \approx B$, we mean that $A\lesssim B$ and $B\lesssim A$.

\section{Preliminaries}\label{sec 2}

For $j\in\mathbb{Z}$, $m\in\mathbb{Z}^{n}$, let $Q_{j,m}:=\prod_{i=1}^{n}[2^{-j}m_{i},2^{-j}(m_{i}+1))$.
For a cube $Q$, $\ell(Q)$ stands for the length of cube $Q$. We
denote by $\mathcal{D}$ the the family of all dyadic cubes in $\mathbb{R}^{n}$,
while $\mathcal{D}_{j}$ is the set of all dyadic cubes with $\ell(Q)=2^{-j},j\in\mathbb{Z}$. Let $x_Q$  be the lower left corner of $Q\in \D$. For $\lambda >0$, let $\lambda Q$ be the cube with the same center of  $Q$ and the edge length $\lambda \ell (Q)$.
For $k\in \mathbb N$, let $k_{\operatorname{pa}} Q$ be the $k$-th dyadic parent of $Q$, which is the dyadic cube in $\D$ satisfying $Q \subset k_{\operatorname{pa}} Q$ and $\ell(   k_{\operatorname{pa}} Q ) = 2^k \ell (Q)$.

\subsection{Matrix weights}
First, we recall some basic concepts and results from the theory of matrix weights.

For any $d,n \in \mathbb N$,  denote by $M_{d,n}  (\mathbb C)$
the set of all $d\times n$ complex-valued matrices. $M_{d,d} (\mathbb C)$ is simply denoted by $M_d (\mathbb C)$. The zero matrix in $M_{d,n}  (\mathbb C) $ (or $M_{d}  (\mathbb C)$) is denoted by $O_{d,n}$ (or $O_d$). Denote by $A^*$ the conjugate transpose of $A \in M_{d,n}  (\mathbb C)$.  A matrix $A \in M_{d}  (\mathbb C) $ is called a {\it Hermitian  matrix } if $A = A^*$ and is called a {\it unitary matrix} if $A^* A = I_d$ where $I_d$ is the identity matrix.
We denote a diagonal matrix by $\operatorname{diag}(\lambda_1, \ldots , \lambda_d) = \operatorname{diag}(\lambda_i).$

For a vector $x\in \mathbb C ^d$, let $|x| = ( \sum_{i=1}^d |x_i|^2 )^{1/2} .$ For $1\le p < \infty$, let $ |x| _p =( \sum_{i=1}^d |x_i|^p )^{1/p}  $. For $ p=\infty$, let $ |x|_\infty = \max(x_1,\ldots, x_d )$. In the finite dimension space, the norms are equivalent. That is, for $ 1\le p<q \le \infty $,
\begin{equation} \label{finite norm equ}
	|x|_q \le |x|_p \le d ^{1/p}  |x|_\infty \le d ^{1/p}  |x|_q.
\end{equation}

For $A \in M_{d}  (\mathbb C)$, let
\begin{equation*}
	\left\| A \right\| := \sup _ {\vec z \in \mathbb C ^d, \left| \vec z \right|  =1 }  \left|  A \vec z  \right|.
\end{equation*}

We say that a matrix  $A  \in M_d(\mathbb{C}) $  is  {\it positive definite}  if, for any  $ \vec x \in \mathbb C ^d \backslash \{\vec 0\} $, $ \vec  z  ^*  A\vec z  >0. $  And a matrix  $A  \in M_d(\mathbb{C}) $  is  called  {\it nonnegative definite} if, for any  $ \vec x \in \mathbb C ^d  $, $ \vec  z  ^*  A\vec z  \ge 0. $

From \cite[Theorem 4.1.4]{HJ13},
any nonnegative definite matrix is always Hermitian. Hence  any nonnegative definite matrix is self-adjoint.

Let $A\in M_d(\mathbb{C})$ be a positive definite matrix
and have eigenvalues $\{\lambda_i\}_{i=1}^d$.
From \cite[Theorem 2.5.6(c)]{HJ13},
there exists a unitary matrix $U\in M_d(\mathbb{C})$ such that
\begin{equation}\label{matrix diag}
	A=U\operatorname{diag}\,(\lambda_1,\ldots,\lambda_d)U^*.
\end{equation}
Moreover, by \cite[Theorem 4.1.8]{HJ13},
we find $\{\lambda_i\}_{i=1}^d\subset(0,\infty) $.
The following definition is based on these conclusions.

\begin{definition}
	Let $A\in M_d(\mathbb{C})$ be a positive definite matrix
	with positive eigenvalues $\{\lambda_i\}_{i=1}^d$.
	For any $\alpha\in\mathbb{R}$, define
	\begin{equation} \label{A alpha}
		A^\alpha:=U\operatorname{diag}\left(\lambda_1^\alpha,\ldots,\lambda_m^\alpha\right)U^*,
	\end{equation}
	where $U$ is the same as in \eqref{matrix diag}.
\end{definition}

\begin{remark}
	From \cite[p.\,408]{HJ94}, we obtain that $A^\alpha$
	is independent of the choices of the order of $\{\lambda_i\}_{i=1}^m$ and $U$,
	and hence $A^\alpha$ is well defined.
\end{remark}

Now, we recall some concepts of matrix wights.
\begin{definition} \label{matrix def}
	A matrix-valued function $W: \rn \to M_d (\mathbb C)$ is called a {\it matrix weight} if $W$ satisfies that
	
	{\rm (i)} for any $x\in \rn$, $W(x)$ is nonnegative definite;
	
	{\rm (ii)} for almost every  $x\in \rn$, $W(x)$ is invertible;
	
	{\rm (iii)} the entries of  $W$ is locally integrable.
\end{definition}

\begin{definition}\label{reduce}
	Let $p\in(0,\infty)$, $W$ be a matrix weight. Suppose that
	$E\subset\mathbb{R}^n$  is a bounded measurable set satisfying $0 < |E| < \infty$. Then
	the matrix $A_E\in M_d(\mathbb{C})$ is called a reducing operator of order $p$ for $W$
	if $A_E$ is positive definite and,
	for any $\vec z\in\mathbb{C}^d$,
	\begin{equation}\label{equ_reduce}
		\left|A_E\vec z\right|
		\approx \left( \frac{1}{|E|} \int_E\left|W^{\frac{1}{p}}(x)\vec z\right|^p\, \d x\right)^{\frac{1}{p}},
	\end{equation}
	where the positive equivalence constants depend only on $d$ and $p$.
\end{definition}

Next we recall the concepts of scalar weight class $A_p$ (see \cite[Definitions 7.1.1, 7.1.3]{Gra1}) and  matrix weight class $\mathcal A_p$ (see \cite{R02} for $1<p<\infty$, \cite{FR04} for $0<p \le 1$).

\begin{definition}
	A weight $\omega$ is a nonnegative  locally integrable function on $\rn$ such that $0<\omega (x) <\infty$ for almost all $x\in \rn$.
	
	A weight $\omega  $ is called an $A_1$ weight if $\mathcal M (\omega) (x) \le c \omega  (x)$ for almost all $x\in \rn$ where $\mathcal M$ is the Hardy-Littlewood maximal operator.
	
	For $1<p<\infty$, a weight $\omega$ is said to be of class $A_p$ if
	\begin{equation*}
		\underset{Q \;\operatorname{cubes\; in}  \rn}{\sup} \left(  \frac{1}{|Q|}\int_Q \omega (x) \d x \right) \left(  \frac{1}{|Q|}\int_Q \omega (x) ^ { -1/ (p-1) } \d x \right)^{p-1} <\infty.
	\end{equation*}
	
	For $1< p < \infty$, a matrix weight $W\in \mathcal A _p(\mathbb{\mathbb{R}}^{n})$
	if and only if
	\[
	\sup_{Q}\frac{1}{|Q|}\int_{Q}\left(\frac{1}{|Q|}\int_{Q}\|W^{1/p}(x)W^{-1/p}(y)\|^{p'}\mathrm{d}y\right)^{p/p'}\mathrm{d}x<\infty,
	\]
	where $p'=p/(p-1)$ is the conjugate index of $p$, and the supremum
	is taken over all cubes $Q\subset\mathbb{R}^{n}$.
	
	For $0<p\le 1$, a matrix weight $W\in \mathcal A _p(\mathbb{\mathbb{R}}^{n})$
	if and only if
	\[
	\underset{Q}{\sup}\;\underset{y\in Q}{\mathop{\mathrm{ess\;sup}}}\frac{1}{|Q|}\int_{Q}\|W^{1/p}(x)W^{-1/p}(y)\|^{p}\mathrm{d}x<\infty.
	\]
	
	We write $\mathcal A _p:=\mathcal A _p(\mathbb{\mathbb{R}}^{n})$ for brevity.
\end{definition}

Given any matrix weight $W$ and $ 0< p <\infty$, there exists (see
e.g., \cite[Proposition 1.2]{G03} for $p>1$ and \cite[p.1237]{FR04}
for $0<p\le1)$ a sequence $\{A_{Q}\}_{Q\in\mathcal{D}}$ of positive
definite $d\times d$ matrices such that
\[
c_{1}|A_{Q}\vec{y}|\le\Big(\frac{1}{|Q|}\int_{Q}|W^{1/p}(x)\vec{y}|^{p}\mathrm{d}x\Big)^{1/p}\le c_{2}|A_{Q}\vec{y}|,
\]
with positive constants $c_{1}$, $c_{2}$ independent of $\vec{y}\in\mathbb{C}^{d}$
and $Q\in\mathcal{D}$. In this case, we call $\{A_{Q}\}_{Q\in\mathcal{D}}$
a sequence of reducing operators of order $p$ for $W$.

\begin{definition}
	\label{def:doubling matrix} A matrix weight $W$ is called a doubling
	matrix weight of order $p>0$ if the scalar measures $w_{\vec{y}}(x)=|W^{1/p}(x)\vec{y}|^{p}$,
	for $\vec{y}\in\mathbb{C}^{d}$, are uniformly doubling: there exists
	$c>0$ such that for all cubes $Q\subset\mathbb{R}^{n}$ and all $\vec{y}\in\mathbb{C}^{d}$,
	\[
	\int_{2Q}w_{\vec{y}}(x)dx\le c\int_{Q}w_{\vec{y}}(x)\mathrm{d}x.
	\]
	If $c=2^{\beta}$ is the smallest constant for which this inequality
	holds, we say that $\beta$ is the doubling exponent of $W$.
	
	From
	\cite[Proposition 2.10]{HS14}, we know that $\beta$ is always not
	less than $n$.
\end{definition}

In \cite{BHYY23}, Bu, Hyt\"onen, Yang and Yuan introduced the $\mathcal A _p$-dimension of matrix weights, which will be used in Theorems  \ref{tm ave BM} and \ref{char Bourgain}.

\begin{definition}
	Let $0<p <\infty$, $\tilde d \in \mathbb R$. A  matrix weight  $W$ has the $\mathcal A _p$-dimension $\tilde d$, denoted by $W \in\mathbb D _{p,\tilde d} ( \mathbb R^n,\mathbb C^n)$, if there exists a positive constant $C$ such that for any cube $Q \subset \mathbb R^n$ and any $i\in \mathbb N  _0$,
	\begin{equation} \label{reducing A_Q}
		\underset{ y\in 2^i  Q}{ \operatorname{ess\;\sup} }  \frac{1}{|Q|} \int_Q \left\| W^{1/p} (x) W^{-1/p} (y) \right\|^p \d x \le C2^{i \tilde d },\quad \operatorname{for} \; 0<p \le 1 ,
	\end{equation}
	or,
	\begin{equation*}
		\frac{1}{|Q|}    \int_Q  \left( \frac{1}{|2^i Q|} \int_{2^i Q} \left\| W^{1/p} (x) W^{-1/p} (y) \right\|^{p'} \d y \right)^{p/p'}  \d x \le C2^{i \tilde d },\quad \operatorname{for} \; 1<p <\infty  ,
	\end{equation*}
	where $ 1/p+ 1/p'= 1$.
\end{definition}
We denote $ \mathbb D _{p,\tilde d}  ( \mathbb R^n,\mathbb C^n) $ simply by $  \mathbb D _{p,\tilde d} $.

The following lemma says that if $W \in  \mathcal A _p $ for $p\in (0,\infty)$, then $W  \in \mathbb D _{p,\tilde d} ( \mathbb R^n,\mathbb C^n)$.
\begin{lemma}[Proposition 2.27, \cite{BHYY23}]
	Let $p \in (0,\infty)$ and $W \in \mathcal A _p$.  Then there exists $\tilde d \in [0,n ) $ such that $W$ has the $\mathcal A _p$-dimension $\tilde d$.
\end{lemma}

\begin{lemma}[Corollary 2.32, \cite{BHYY23}] \label{impro AQ AR}
	Let $0<p<\infty$, let $W \in \mathcal A _p$ with the $\mathcal A _p$-dimension $\tilde d  \in [0,n)$, and let $\{ A_Q\}_ { \operatorname{cube} Q } $ be the reducing operator of order $p$  for $W$.
	{\rm (i)} If $1<p<\infty$, let $\widetilde W := W ^{-1/(p-1)}$ (which belongs to $\mathcal A _{p'}$) with the $\mathcal A _{p'}$-dimension $\tilde {\tilde d}$.
	Then there exists a positive constant $C$ such that, for any cubes $Q$ and $R$ of $\rn$,
	\begin{equation*}
		\| A_Q A_R ^{-1} \|\le C \max \left( \left[ \frac{\ell(R)}{\ell(Q)}\right]^{\tilde d /p}, \left[ \frac{\ell(Q)}{\ell(R)}\right]^{\tilde {\tilde d} /p'}   \right)  \left[ 1+ \frac{|x_Q - x_R| }{  \max(\ell(Q) , \ell (R)  ) }\right] ^{ \tilde d /p + \tilde {\tilde d} /p' }.
	\end{equation*}
	{\rm (ii)} If $0<p\le 1$,
	then there exists a positive constant $C$ such that, for any cubes $Q$ and $R$ of $\rn$,
	\begin{equation*}
		\| A_Q A_R ^{-1} \|\le C \max \left(  \left[ \frac{\ell(R)}{\ell(Q)}\right]^{\tilde d /p}, 1   \right)  \left[ 1+ \frac{|x_Q - x_R| }{  \max(\ell(Q) , \ell (R)  ) }\right] ^{ \tilde d /p }.
	\end{equation*}
	
\end{lemma}

\subsection{Matrix weighted Bourgain-Morrey spaces}

\begin{definition}
	Let $0 < p <\infty$, $d\in \mathbb N$ and $W:\rn \to M_d (\mathbb C)$ be a matrix weight. The space $L^p  _{\operatorname{loc}}(W)$  collects all measurable functions $\vec f : \rn \to \mathbb C ^d$ such that for each compact set $K$,
	\begin{equation*}
		\| \vec f \chi_K \|_{L^p (W) }  = \left(   \int_K | W^{1/p} (x) \vec f (x) |^p \d x \right)^{1/p}  <\infty.
	\end{equation*}
	Let $\Omega \subset \rn$ be an open set.
	The space $L^p (W,\Omega)$  collects all measurable functions $\vec f : \Omega \to \mathbb C ^d$ such that
	\begin{equation*}
		\| \vec f  \|_{L^p (W,\Omega) }  = \left(   \int_\Omega | W^{1/p} (x)\vec f (x) |^p \d x \right)^{1/p}  <\infty.
	\end{equation*}
\end{definition}

\begin{definition}
	Let $\D  = \{ Q_{j,k} \}_{ j\in \mathbb Z , k\in \mathbb Z ^n} $ be the standard dyadic system.
	Let $0<p<  t <r <\infty$ or $0< p \le t< r =\infty $. Let  $W :\rn \to M_d (\mathbb C)$  be a matrix weight. Define $ M_{p}^{t,r} (W) $ as the set of all $\vec f \in L _{\rm loc}^p (W)$ such that
	\begin{equation*}
		\| \vec f \|_{  M_{p}^{t,r} (W)    } := \left\|     \left\{    W (Q_{j,k} )^{ 1/t-1/p }   \left(  \int_{Q_{j,k}  }  | W^{1/p} (y)  \vec f (y) |^p \d y     \right) ^{1/p}    \right\}    _{ j\in \mathbb Z , k\in \mathbb Z ^n }   \right\|   _{ \ell^r } <\infty,
	\end{equation*}
	where $  W (Q_{j,k} ) = \int_{Q_{j,k} } \| W (y) \| \d y $.
	
	Let  $\{ A_Q\}_ {Q \in \D} $ be the reducing operator of order $p$  for $W$.  Define $M_{p}^{t,r} (  \{A_Q\}  )$ as the set of all $\vec f \in L _{\rm loc}^p (W)$ such that
	\begin{align*}
		\| \vec f \|_{  M_{p}^{t,r} ( \{A_Q\}  )    } :=  \left\|     \left\{    \left( \|  A_{Q_{j,k} } \|^p |Q_{j,k}  | \right)^{ 1/t-1/p }   \left(  \int_{Q_{j,k}  }  | W^{1/p} (x) \vec f (y) |^p \d y     \right) ^{1/p}    \right\}    _{ j\in \mathbb Z , k\in \mathbb Z ^n }   \right\|   _{ \ell^r }  <\infty.
	\end{align*}
\end{definition}
\begin{remark}
	By \cite[Lemma 2.11]{BHYY23} with $M = I_m$, we conclude that, for any cube $Q \subset \rn$,
	\begin{equation*}
		\|A_Q\|^p \approx \frac{1}{|Q|}  \int_Q  \left\| W^{1/p} (x) \right\|^p \d x=  \frac{1}{|Q|} \int_Q  \left\| W (x) \right\| \d x= \frac{1}{|Q|} W (Q).
	\end{equation*}
	Thus $M_{p}^{t,r} (W)$ is same with $ M_{p}^{t,r} (  \{A_Q\}  )$ in meaning of the  equivalent quasi-norms.
	
	If $ d = 1, W \equiv 1 $, $ M_{p}^{t,r} (W) $ is the classical Bourgain-Morrey space $M_p^{t,r}$ in \cite{HNSH23}.
	We define  the scalar weighted  Bourgain-Morrey space $M_{p}^{t,r} (\omega)   $ by
	\begin{equation} \label{scalar wei BM}
		\| f \|_{ M_p^{t,r} (\omega) } :=   \left\|     \left\{  (    \omega ( Q_{j,k} )  )^{ 1/t-1/p }   \left(  \int_{Q_{j,k}  }  |  f (y) |^p \omega (y) \d y     \right) ^{1/p}    \right\}    _{ j\in \mathbb Z , k\in \mathbb Z ^n }   \right\|   _{ \ell^r } <\infty.
	\end{equation}
	That is  $M_{p}^{t,r} (\omega)   $ is the case $d=1 $ of matrix weighted Bourgain-Morrey space $M_{p}^{t,r} (W)  $.
	
\end{remark}

The following lemma is proved in \cite{MRL12,SW10}.

\begin{lemma} \label{Banach Lp Omega}
	Let $1\le p <\infty$ and $W:\rn \to M_d (\mathbb C)$ be a matrix weight.
	Let $\Omega \subset \rn$ be an open set.
	Then matrix weighted Lebesgue space $L^p(W, \Omega)$ is a  Banach space.
\end{lemma}

\begin{theorem}
	Let $1\le p <\infty$ and $W:\rn \to M_d (\mathbb C)$ be a matrix weight. Then the space $L^p  _{\operatorname{loc}}(W)$ is complete.
\end{theorem}

\begin{proof}
	Let $\{\vec f _k \} _{ k =1 } ^\infty$ be a Cauchy sequence in $L^p  _{\operatorname{loc}}(W)$.
	That is, for any compact set $K$, any $\epsilon >0$, there exists $N >0$ such that if $j,k >N$, $	\|  ( \vec f_j - \vec f _k ) \chi_K \| _ {L^p (W)}  < \epsilon .$
	Then for any  compact set $K$, $ \{\vec f _k  \chi_K \} _{ k =1 } ^\infty $ is a Cauchy sequence in $L^p (W)$. Since $L^p (W)$ is complete (Lemma \ref{Banach Lp Omega}), there exists $\vec f \chi_K$ in $L^p (W)$ such that $\vec f _k  \chi_K \to \vec f \chi_K$. By the Fatou lemma, we have
	\begin{equation*}
		\| \vec f \chi_K \| _ {L^p (W)} \le \liminf_{k\to \infty} 	\| \vec f _k \chi_K \| _ {L^p (W)}<\infty.
	\end{equation*}
	By the dominated convergence theorem, for $j,k>N$,  we obtain
	\begin{equation*}
		\| (\vec f   - \vec f_k) \chi_K \| _ {L^p (W)} = \lim_{j\to \infty}  \|  (\vec f_j   - \vec f_k) \chi_K \| _ {L^p (W)} <\epsilon.
	\end{equation*}
	Hence $ \vec f_k \to \vec f $ in $L^p  _{\operatorname{loc}}(W)$ and the space $L^p  _{\operatorname{loc}}(W)$ is complete.
\end{proof}

In what follows, the symbol $\hookrightarrow $ always stands for continuous embedding.

\begin{proposition} \label{basic embedding}
	Let  $W: \rn \to M_d (\mathbb C)  $ be a matrix weight.
	
	{\rm (i)} If $ 0<p < t < r_1 < r_2 \le \infty $, then
	\begin{equation*}
		M_{p}^{t,r_1} (W) \hookrightarrow  M_{p}^{t,r_2} (W) .
	\end{equation*}
	
	{\rm (ii)} If $ 0 <p_1  <p_2 <t <r<\infty $ or  $ 0 <p_1 <p_2 \le t < r =\infty $, then
	\begin{equation*}
		M_{p_2}^{t,r} (W) \hookrightarrow  M_{p_1}^{t,r} (W) .
	\end{equation*}
	
	{\rm (iii)} If  $ 0 <p \le t < r =\infty $, then
	\begin{equation*}
		L^ t (W) \hookrightarrow	M_{p}^{t,r} (W) \hookrightarrow L_{\operatorname{loc}}^p (W) .
	\end{equation*}
	
\end{proposition}
\begin{proof}
	(i) It comes from $ \ell^{r_1}  \hookrightarrow \ell ^{r_2}$ since $0<r_1 <r_2 \le \infty $.
	
	(ii) It comes from the H\"older  inequality. Indeed, for each $Q \in \D$,
	\begin{align*}
		&	\left(  \frac{1}{W (Q)} \int_Q  |W^{1/p_1} (y) \vec  f (y) |^{p_1} \d y  \right)^{1/p_1} \\
		& =\left(  \frac{1}{W (Q)} \int_Q  |W^{1/p_1  - 1/p_2} (y)  W^{1/p_2} (y)\vec  f (y) |^{p_1} \d y  \right)^{1/p_1} \\
		& \le \left(  \frac{1}{W (Q)} \int_Q  \|W^{1/p_1  - 1/p_2} (y)\| ^{p_1}   | W^{1/p_2} (y)\vec  f (y) |^{p_1} \d y  \right)^{1/p_1} \\
		& \le \left(  \frac{1}{W (Q)} \int_Q  \|W^{1/p_1  - 1/p_2} (y)\| ^{p_1 p_2/ (p_2 - p_1 )} \d y  \right)^{1/p_1-1/p_2}
		\left( \frac{1}{W (Q)} \int_Q   | W^{1/p_2} (y)\vec  f (y) |^{p_2} \d y  \right)^{1/p_2} \\
		& = \left(  \frac{1}{W (Q)} \int_Q  \|W (y)\|  \d y  \right)^{1/p_1-1/p_2}  \left( \frac{1}{W (Q)} \int_Q   | W^{1/p_2} (y)\vec  f (y) |^{p_2} \d y  \right)^{1/p_2} \\
		&=\left( \frac{1}{W (Q)} \int_Q   | W^{1/p_2} (y)\vec  f (y) |^{p_2} \d y  \right)^{1/p_2}.
	\end{align*}
	Hence,
	\begin{align*}
		\| \vec f \| _{M_{p_1}^{t,r} (W)  }\le  \| \vec f \| _{M_{p_2}^{t,r} (W)  }.
	\end{align*}
	Thus we prove (ii).
	
	(iii) The first embedding comes from the fact that
	\begin{equation*}
		L^ t (W) =	M_{t}^{t,\infty} (W)    \hookrightarrow  M_{p}^{t,\infty} (W)   .
	\end{equation*}
	For any compact set $K\subset \rn$, there exist at most $2^n$ dyadic cubes $Q_j$ such that $ K \subset \bigcup_{j=1}^{2^n} Q_j $. Hence
	\begin{align*}
		\| \vec f \chi_K \|_{L^p (W)} \le \sum_{j=1}^{2^n}  \frac{|Q_j|^{1/t-1/p} }{|Q_j|^{1/t-1/p}} \| \vec f \chi_{Q_j} \|_{L^p (W)}
		\le \| \vec f \|_{M _P^{t,\infty} (W) } \sum_{j=1}^{2^n}  \frac{1 }{|Q_j|^{1/t-1/p}} <\infty.
	\end{align*}
	This shows that 	$M_{p}^{t,r} (W) \hookrightarrow L_{\operatorname{loc}}^p (W) .$
\end{proof}

\begin{proposition}
	The matrix weighted Bourgain-Morrey space $M_{p}^{t,r} (W)   $ is a Banach space when $1\le p <t <r <\infty$  or $ 1\le p \le  t <r=\infty$.
	
\end{proposition}
\begin{proof}
	The argument are standard, see, for example, \cite[Theorem 2.4]{S18}.
	We only show that  $M_{p}^{t,r} (W) $ is complete since others are simple.

	Let $\{ \vec f _j \}_{j=1}^\infty$  be a Cauchy sequence in $M_{p}^{t,r} (W) $. By Proposition \ref{basic embedding}, we have that $\{ \vec f _j \}_{j=1}^\infty$  is also a Cauchy sequence in $L^p_{ \operatorname{loc}  } (W) $. Since $L^p_{ \operatorname{loc}  } (W) $  is complete, there exists a vector function $\vec f$ in $L^p_{ \operatorname{loc}  } (W) $ such that $  \vec f _j \to \vec f $.
	By the Fatou lemma,
	\begin{align*}
		\| \vec f \|_{ M_{p}^{t,r} (W) }  = 	\|  \underset{j\to \infty }{\lim \inf} \vec f_j \|_{ M_{p}^{t,r} (W) } \le   \underset{j\to \infty }{\lim \inf} 	\|   \vec f_j \|_{ M_{p}^{t,r} (W) }  <\infty.
	\end{align*}
	Thus $ \vec f \in  M_{p}^{t,r} (W) $. Therefore,
	\begin{align*}
		\underset{j\to \infty }{\lim \sup} 	\| \vec f  - \vec f_j \|_{ M_{p}^{t,r} (W) }  \le  \underset{j\to \infty }{\lim \sup}  \left(  \underset{k \to \infty }{\lim \inf}  	\| \vec f_k  - \vec f_j \|_{ M_{p}^{t,r} (W) }     \right) =0.
	\end{align*}
	Thus, the sequence  $\{ \vec f _j \}_{j=1}^\infty$ is convergent to $\vec f$ in $M_{p}^{t,r} (W)$.
\end{proof}

Recall that a matrix weight $W$ is almost everywhere invertible. Hereafter, we define $ w (x) :=\| W ^{-1} (x) \|^{-1}$ when a matrix weight $W$   is  invertible at $x \in \rn$.
\begin{lemma}[Proposition 3.2,\cite{CMR16}]  \label{ellipticity}
	Let $ 1\le p <\infty$.
	For a matrix weight $W:  \rn \to M_d (\mathbb C)$, we have $ 0 < w (x) \le  \|  W (x) \| <\infty$  for a.e. $x \in \rn$. Furthermore, $W$ satisfies a two weight, degenerate ellipticity 	condition: for $ \xi  \in \mathbb R ^d$,
	\begin{equation} \label{deg ell con}
		w (x) |\xi |^p \le | W^{1/p} (x) \xi |^p \le \| W(x)\| |\xi |^p .
	\end{equation}
\end{lemma}
\begin{remark}
	In \cite[Proposition 3.2]{CMR16}, it is assumed that $\xi \in \mathbb R^d$ but it also works for  $\xi \in \mathbb C^d$.
\end{remark}

\begin{proposition}
	Let $1\le p <t <r <\infty$  or $ 1\le p \le  t <r=\infty$. Let $W$  be a matrix weight and $ \vec f, \vec g \in M_{p}^{t,r} (W)$. Then $ \| \vec f - \vec g \|_{ M_{p}^{t,r} (W)  }  = 0$ if and only if $ \vec f (x) =\vec g  (x)$ a.e..
\end{proposition}
\begin{proof}
	Clearly, if   $ \vec f (x) =\vec g  (x)$ a.e., then  $ \| \vec f - \vec g \|_{ M_{p}^{t,r} (W)  }  = 0$.
	Then we apply Lemma \ref{ellipticity} to prove the converse.
	By the degenerate ellipticity 	condition (\ref{deg ell con}), we have
	\begin{align*}
		0 = \| \vec f - \vec g \|_{ M_{p}^{t,r} (W)  } \ge  \left\|     \left\{    W (Q_{j,k} )^{ 1/t-1/p }   \left(  \int_{Q_{j,k}  }  |   \vec f (y) -\vec g (y)|^p  w (y) \d y     \right) ^{1/p}    \right\}    _{ j\in \mathbb Z , k\in \mathbb Z ^n }   \right\|   _{ \ell^r }.
	\end{align*}
	Since $w(y) >0 $ a.e., it follows that $\vec f (y)-\vec g(y) = 0$ a.e..
\end{proof}

\section{Sufficient conditions for precompact sets}\label{KRct}

\begin{definition}
	{\rm(i)} Suppose that  $ (X,\rho )$ is  a metric space. Let $A$ be a set of $X$  and let $\epsilon >0$. A set $ E \subset X$ is called an $ \epsilon$-net for $A$ if every point $a\in A$, there exists a point  $ e\in E $ such that $ \rho(a,e) <\epsilon $.
	
	{\rm (ii)} A set $A$ is called totally bounded if, for each $\epsilon>0$, it possesses a finite $\epsilon$-net.
	
	{\rm (iii)} A subset in a topological space is precompact if its closure is compact.
\end{definition}
It is well known that  in a complete metric space, a set is precompact if and only if it is totally
bounded.

\begin{theorem} \label{K-R 1}
	Let $1 \le  p < t < r <\infty$ or $ 1 \le p \le t < r= \infty $. Let $W: \rn \to  M_d (\mathbb C) $ be a matrix weight.
	A subset $\mathcal F \subset M_{p}^{t,r} (W) $ is totally bounded if the following conditions are valid:

	{\rm  (i)} $\mathcal F$ uniformly vanishes at infinity, that is,
	\begin{equation*}
		\lim_{R\to \infty} \sup _{\vec f \in  \mathcal F} 	\| \vec f \chi_{B^c (0,R)}  \|_{  M_{p}^{t,r} (W)    }  =0;
	\end{equation*}
	
	{\rm  (ii)} $\mathcal F$ is equicontinuous, that is,
	\begin{equation}\label{trans con}
		\lim _{b \to 0} \sup _{\vec f \in  \mathcal F} \sup_{y\in B(0,b)} 	\| \vec f  - \tau_y \vec  f  \|_{ M_{p}^{t,r} (W)  } =0.
	\end{equation}
	Here and what follows, $\tau _y$ denotes the translation operator: $\tau_y \vec f (x) := \vec f (x-y)$.
\end{theorem}

\begin{remark}\label{F bounded}
	In this Remark, we will prove that	the conditions (i) and (ii)  in Theorem \ref{K-R 1} together imply that the set $\mathcal F$ is bounded. Indeed,	choose $b>0$ such that for all $h\in \rn, |h| \le b$, all $\vec f\in \mathcal F$,
	\begin{equation*}
		\| \vec f - \tau_y \vec f \|_{  M ^{t,r}_{p} (W) }  \le 1.
	\end{equation*}
	Choose $R>0$  such that for all $\vec f\in \mathcal F$,
	\begin{equation*}
		\| \vec f\chi_{ B^c (0,R) } \|_{  M ^{t,r}_{p}  (W) }  \le 1.
	\end{equation*}
	Fix $h$ with $|h| =b$.
	Then for all $\vec f\in \mathcal F$, $x\in \rn$,  we have
	\begin{align*}
		\| \vec f \chi_{ B(x,R) } \|_{ M ^{t,r}_{p}(W)   } & \le  	\| ( \vec f- \tau_h \vec f ) \chi_{ B(x,R) } \|_{ M ^{t,r}_{p} (W)  }  + \|  \tau_h \vec f  \chi_{ B(x,R) } \|_{ M ^{t,r}_{p} (W)  }  \\
		& =	\| ( \vec f- \tau_h \vec f ) \chi_{ B(x,R) } \|_{ M ^{t,r}_{p} (W)  }  + \|  \vec f  \chi_{ B(x+ h,R) } \|_{ M ^{t,r}_{p} (W)  }  \\
		&\le 1 + \|  \vec f  \chi_{ B(x+ h,R) } \|_{ M ^{t,r}_{p} (W)  }  .
	\end{align*}
	Hence, by induction,
	\begin{equation*}
		\| \vec f \chi_{ B(0,R) } \|_{ M ^{t,r}_{p}(W)   } \le N + \|  \vec f  \chi_{ B(N h,R) } \|_{ M ^{t,r}_{p} (W)  }  .
	\end{equation*}
	Now choose $N \ge 1$ such that $ N b = N |h| > 2R $. Then $ B(N h , R)  \subset B^c (0,R) $.
	Hence
	\begin{align*}
		\| \vec f \|_{ M ^{t,r}_{p} (W) } \le \| \vec f \chi_{ B(0,R) } \|_{ M ^{t,r}_{p}(W)   }  + \| \vec f \chi_{ B^c(0,R) } \|_{ M ^{t,r}_{p}(W)   }
		\le N + \|  \vec f  \chi_{ B(N h,R) } \|_{ M ^{t,r}_{p} (W)  }  + \| \vec f \chi_{ B^c(0,R) } \|_{ M ^{t,r}_{p}(W)   }  \le N+2.
	\end{align*}
	This proves that $\mathcal F$ is bounded.
\end{remark}
Now we begin to show Theorem \ref{K-R 1}.
\begin{proof}
	Assume that $\mathcal F \subset  M_{p}^{t,r} (W) $ satisfies (i) and (ii). Given $\epsilon >0$ small enough, to prove the total boundedness of $\mathcal F$, it suffices to find a finite  $\epsilon$-net of $\mathcal F$. Denote by $R_i := [-2^i , 2^i )^n$ for $i\in\mathbb Z$. Then from condition (i), there exist a positive integer $m$ large enough such that
	\begin{equation} \label{epsilon 1}
		\sup _{\vec f \in  \mathcal F} 	\| \vec f - \vec f  \chi_{R_m}  \|_{  M_{p}^{t,r} (W)   }  <\epsilon.
	\end{equation}
	Moreover, by condition (ii), there exists a negative integer $a$ such that
	\begin{equation} \label{condi 2}
		\sup _{\vec f \in  \mathcal F} \sup_{y\in R_a} 	\| \vec f  - \tau_y \vec  f  \|_{  M_{p}^{t,r} (W)   }  <\epsilon.
	\end{equation}
	There exists a sequence $\{ Q_j \} _{j=1}^N$ of disjoint cubes in $\D _{-a}$ such that $R_m = \bigcup_{i=1}^N  Q_j$, where $N=2^{  (m+1-a ) n}$.
	For any $\vec f \in \mathcal F$ and $x\in \rn$, let
	\begin{equation*}
		\Phi (\vec{ f}) (x):= \begin{cases}
			\vec f_{Q_j} :=\frac{1}{|Q_j|} \int_{Q_j} \vec f (y) \d y, &  x\in Q_j, j = 1,2,\ldots, N, \\
			\vec 0,  & \operatorname{otherwise}.
		\end{cases}
	\end{equation*}
	Then for each fixed  $x\in \rn$, we have
	\begin{align}
		\label{f Qj}
		\nonumber
		\left|	 W^{1/p}(x) \left( \vec f (x) - 	\vec f_{Q_j} \right)    \right| \chi_{Q_j} (x) & = \left| \frac{1}{|Q_j|} \int_{Q_j} W^{1/p}(x) \left( \vec f (x) - 	\vec f (y) \right )  \d y   \right|   \chi_{Q_j} (x) \\
		& \le \frac{1}{|Q_j|} \int_{Q_j}  \left|   W^{1/p}(x) \left ( \vec f (x) - 	\vec f (y)  \right )     \right|  \d y  \chi_{Q_j} (x).
	\end{align}
	We split 	$ \|  ( \vec f - \Phi( \vec f)  ) \chi_{R_m} \|_{ M_{p}^{t,r} (W) } $ into three parts:
	\begin{align*}
		&	\left\{	\sum_{k\in \mathbb Z} \sum_{Q \in \mathcal D_k} W(Q)^{r/t - r/p} \left( \int_Q 	\left|	 W^{1/p}(x) \left ( \vec f (x) \chi_{R_m}- 	\Phi (\vec f) (x) \right)    \right| ^p \d x \right) ^{r/p}  \right\}^{1/r} \\
		&\le 	\left\{	\sum_{k > -a} \sum_{Q \in \mathcal D_k} W(Q)^{r/t - r/p} \left( \int_Q 	\left|	 W^{1/p}(x) \left ( \vec f (x)  \chi_{R_m} - 	\vec f_{Q_j} \right )    \right| ^p \d x \right) ^{r/p}  \right\}^{1/r} \\
		&	\quad +	\left\{	\sum_{k = -m }^{-a } \sum_{Q \in \mathcal D_k} W(Q)^{r/t - r/p} \left( \int_Q 	\left|	 W^{1/p}(x) \left ( \vec f (x)  \chi_{R_m} - \Phi (\vec f) (x) \right )    \right| ^p \d x \right) ^{r/p}  \right\}^{1/r} \\
		&	\quad +	\left\{	\sum_{k <-m} \sum_{Q \in \mathcal D_k} W(Q)^{r/t - r/p} \left( \int_Q 	\left|	 W^{1/p}(x) \left ( \vec f (x)  \chi_{R_m} - \Phi (\vec f) (x)\right )    \right| ^p \d x \right) ^{r/p}  \right\}^{1/r} \\
		&  =: S_1 +S_2+S_3.
	\end{align*}
	We first estimate $S_1$. By (\ref{f Qj}), Jensen's inequality ($p\ge 1$), the Fubini theorem, we have
	
	\begin{align}\label{eq 1}
		\nonumber
		S_1  & \le 	\left\{	\sum_{k > -a} \sum_{Q \in \mathcal D_k, Q \subset R_m} W(Q)^{r/t - r/p}
		\left( \int_Q 	\left| \frac{1}{|Q_j|}	\int_{Q_j}  \left|   W^{1/p}(x) \left ( \vec f (x) - 	\vec f (y)  \right)     \right|  \d y  \chi_{Q_j} (x)   \right| ^p \d x \right) ^{r/p}  \right\}^{1/r} \\
		\nonumber
		& \le \left\{	\sum_{k > -a} \sum_{Q \in \mathcal D_k, Q \subset R_m} W(Q)^{r/t - r/p}
		\left(  \frac{1}{|Q_j|} \int_{Q_j} 	 \int_{Q}  \left|   W^{1/p}(x) \left ( \vec f (x) - 	\vec f (y)  \right )     \right| ^p  \d x   \d y \right) ^{r/p}  \right\}^{1/r}  \\
		\nonumber
		& = \left\{	\sum_{k > -a} \sum_{Q \in \mathcal D_k, Q \subset R_m} W(Q)^{r/t - r/p}
		\left(  \frac{1}{|Q_j|} 	 \int_{Q}   \int_{Q_j}  \left|   W^{1/p}(x) \left( \vec f (x) - 	\vec f (y)  \right )     \right| ^p     \d y \d x \right) ^{r/p}  \right\}^{1/r}  \\
		& =2^{  -an/p } \left\{	\sum_{k > -a} \sum_{Q \in \mathcal D_k, Q \subset R_m} W(Q)^{r/t - r/p}
		\left( 	 \int_{Q}   \int_{ x - Q_j}  \left|   W^{1/p}(x) \left ( \vec f (x) - 	\vec f (x- y)  \right )     \right| ^p     \d y \d x \right) ^{r/p}  \right\}^{1/r}
	\end{align}
	where $x-Q_j :=\{ x-y:y\in Q_j \}$. Since $x \in Q \subset Q_j$ for some  $j  \in\{1,2,\ldots, N\}$, we have $x-Q_j \subset R_a$.
	Hence by (\ref{condi 2}), we obtain
	\begin{align*}
		(\ref{eq 1}) \le 2^{  -an/p }| Q_j | ^{1/p} \sup_{y\in R_a}  \left\{	\sum_{k > -a} \sum_{Q \in \mathcal D_k, Q \subset R_m} W(Q)^{r/t - r/p}
		\left( 	 \int_{Q}     \left|   W^{1/p}(x) ( \vec f (x) - 	\vec f (x- y)  )     \right| ^p   \d x \right) ^{r/p}  \right\}^{1/r}
		\le  \epsilon.
	\end{align*}
	As for $S_2$, for each $Q\in \D_k$ where $ k=-m , -m+1, \ldots, -a$, there are $2^{(-k-a)n}$ cubes $Q_j \in \D _{-a}$ such that $ \cup  Q_{j }  =Q$.
	Denote by $Q_{j_ {\ell}  }$, $\ell = 1,2,\ldots, 2^{(-k-a)n}$ these cubes. Then by (\ref{f Qj}), Jensen's inequality ($p\ge 1$), and the Fubini theorem,  we obtain
	\begin{align}\label{S2 1}
		\nonumber
		S_2 & = 	\left\{	\sum_{k = -m}^{-a} \sum_{Q \in \mathcal D_k, Q \subset R_m} W(Q)^{r/t - r/p}
		\left(
		\sum_{\ell=1}^{  2^{(-k-a)n} }
		\int_{Q_{j_\ell}}\left|	 W^{1/p}(x) \left ( \vec f (x)  \chi_{R_m} - 	\vec f_{Q_ {j_\ell }  } \right )    \right| ^p \d x \right) ^{r/p}  \right\}^{1/r}  \\
		\nonumber
		& \le 	\left\{	\sum_{k = -m}^{-a} \sum_{Q \in \mathcal D_k, Q \subset R_m} W(Q)^{r/t - r/p}
		\left(
		\sum_{\ell=1}^{  2^{(-k-a)n} }
		\int_{Q_{j_\ell}} \left|	\frac{1}{|Q_{j_\ell}|} \int_{Q_ {j_\ell } }  \left|   W^{1/p}(x) \left ( \vec f (x) - 	\vec f (y)  \right)     \right|   \d y \right|^p \d x \right) ^{r/p}  \right\}^{1/r}  \\
		\nonumber
		& \le 	\left\{	\sum_{k = -m}^{-a} \sum_{Q \in \mathcal D_k, Q \subset R_m} W(Q)^{r/t - r/p}
		\left(
		\sum_{\ell=1}^{  2^{(-k-a)n} }
		\frac{1}{|Q_{j_\ell}|} 	\int_{Q_ {j_\ell} }
		\int_{Q_{j_\ell}} \left|   W^{1/p}(x) \left ( \vec f (x) - 	\vec f (y)  \right)     \right| ^p  \d x \d y  \right) ^{r/p}  \right\}^{1/r}  \\
		\nonumber
		& =2^{-an/p}	\left\{	\sum_{k = -m}^{-a} \sum_{Q \in \mathcal D_k, Q \subset R_m} W(Q)^{r/t - r/p}
		\left(
		\sum_{\ell=1}^{  2^{(-k-a)n} }
		\int_{Q_{j_\ell }}
		\int_{Q_{j_\ell}} \left|   W^{1/p}(x) \left ( \vec f (x) - 	\vec f (y)  \right)     \right| ^p \d y \d x   \right) ^{r/p}  \right\}^{1/r}  \\
		& =2^{-an/p}	\left\{	\sum_{k = -m}^{-a} \sum_{Q \in \mathcal D_k, Q \subset R_m} W(Q)^{r/t - r/p}
		\left(
		\sum_{\ell=1}^{  2^{(-k-a)n} }
		\int_{Q_{j_\ell }}
		\int_{x  - Q_{j_\ell}} \left|   W^{1/p}(x) \left ( \vec f (x) - 	\vec f (x-y)  \right)     \right| ^p \d y \d x   \right) ^{r/p}  \right\}^{1/r}
	\end{align}
	where $x-Q_{ j _\ell } :=\{ x-y:y\in Q_{j_\ell } \}$. Note that $x-Q_{ j _\ell } \subset R_a $ when $x\in Q_{ j _\ell } $. By (\ref{condi 2}), we have
	\begin{align*}
		(\ref{S2 1})  \le
		2^{-an/p}  2^{an/p}  \sup_{y\in R_a}	\left\{	\sum_{k = -m}^{-a} \sum_{Q \in \mathcal D_k, Q \subset R_m} W(Q)^{r/t - r/p}
		\left(
		\sum_{\ell=1}^{  2^{(-k-a)n} }
		\int_{Q_{j_\ell }}
		\left|   W^{1/p}(x) \left ( \vec f (x) - 	\vec f (x-y)  \right)     \right| ^p \d x   \right) ^{r/p}  \right\}^{1/r}
		\le \epsilon.
	\end{align*}
	As for $S_3$, for each set $\D_k$, there are only $2^n$ cubes $Q \in \D_k$ such that $Q \cap R_m \neq \emptyset $. We denote by $Q_{ \ell }$,  $ \ell =1,2,\ldots , 2^n $ these cubes.
	And for each cube $ Q_\ell $, there are   $2^{ (m-a)n }$ cubes $Q_j \in \mathcal D _{-a}  $ such that $ \cup_j Q_j = Q_\ell $.
	Then
	\begin{align*}
		S_3 &= \left\{	\sum_{k <-m}
		\sum_{\ell=1}^{ 2^n }
		W ( Q _\ell ) ^{r/t - r/p} \left( \int_{Q_\ell}	\left|	 W^{1/p}(x) \left ( \vec f (x)  \chi_{R_m} - 	\Phi (\vec f) (x) \right )    \right| ^p \d x \right) ^{r/p}  \right\}^{1/r}
		\\
		&= \left\{	\sum_{k <-m}
		\sum_{\ell=1}^{ 2^n }
		W ( Q _\ell ) ^{r/t - r/p} \left(  \sum_{j=1}^{ 2^{  (m-a )n} } \int_{Q_j }	\left|	 W^{1/p}(x) \left ( \vec f (x)  \chi_{R_m} - 	\vec f_{Q_j} \right )    \right| ^p \d x \right) ^{r/p}  \right\}^{1/r} .
	\end{align*}
	Similarly as $ S_2 $,  from (\ref{f Qj}), Jensen's inequality ($p\ge 1$), the Fubini theorem, and (\ref{condi 2}), we get that $ S_3 \le \epsilon $.
	Together with the estimates of $S_1$,  $S_2$, $S_3$, we obtain
	\begin{equation} \label{epsilon 2}
		\left\{	\sum_{k\in \mathbb Z} \sum_{Q \in \mathcal D_k} W(Q)^{r/t - r/p} \left( \int_Q 	\left|	 W^{1/p}(x) ( \vec f (x) \chi_{R_m}- 	\Phi (\vec f) (x)  )    \right| ^p \d x \right) ^{r/p}  \right\}^{1/r} \le 3 \epsilon.
	\end{equation}
	Note that
	\begin{align*}
		\| \vec f - \Phi( \vec f) \|_{ M_{p}^{t,r} (W) }  \le 	\|  ( \vec f - \Phi( \vec f)  ) \chi_{R_m} \|_{ M_{p}^{t,r} (W) }  +	\|  ( \vec f - \Phi( \vec f)  ) \chi_{R_m^c} \|_{ M_{p}^{t,r} (W) }
		= 	\|  ( \vec f - \Phi( \vec f)  ) \chi_{R_m} \|_{ M_{p}^{t,r} (W) }  +	\|   \vec f  \chi_{R_m^c} \|_{ M_{p}^{t,r} (W) } .
	\end{align*}
	Thus via (\ref{epsilon 1}) and (\ref{epsilon 2}), we have
	\begin{equation} \label{epsilon Phi}
		\sup_{\vec f \in \mathcal F} 		\| \vec f - \Phi( \vec f) \|_{ M_{p}^{t,r} (W) } \le 4 \epsilon.
	\end{equation}
	From (\ref{epsilon Phi}), it suffices to show that $\Phi( \mathcal F)$ is totally bounded in $M_{p}^{t,r} (W)$.
	
	From Remark \ref{F bounded}, we have
	\begin{equation*}
		\sup_{\vec f \in \mathcal F} \| \Phi (\vec f) \|_{M_{p}^{t,r} (W)  } \le	\sup_{\vec f \in \mathcal F} \| \Phi (\vec f) - \vec  f \|_{M_{p}^{t,r} (W)  } + 	\sup_{\vec f \in \mathcal F} \| \vec f \|_{M_{p}^{t,r} (W)  } <\infty.
	\end{equation*}
	Thus, for any $\vec f \in \mathcal F$,
	\begin{equation}
		|W^{1/p} (x) \Phi (\vec f) (x) |  <\infty  \quad {\rm a.e.} \quad  x\in \rn.
	\end{equation}
	Since $W$ is a matrix weight, by (i), (ii) of Definition \ref{matrix def}, we have for almost everywhere $x\in \rn$, $W(x)$ is a positive definite matrix.
	Therefore, we obtain
	\begin{equation}
		| \Phi (\vec f) (x) |  <\infty  \quad {\rm a.e.} \quad  x\in \rn.
	\end{equation}
	which implies $|\vec f _{Q_j} | <\infty$, $j=1,2,\ldots , N$. From this and the entries of $W$ is locally integrable, we see that $\Phi$ is a map from $\mathcal F$ to $\mathcal B$, a finite dimensional Banach subspace of $M_{p}^{t,r} (W)   $.  Note that $\Phi  (\mathcal F) \subset \mathcal B $ is bounded, and hence is totally bounded. The proof of Theorem \ref{K-R 1} is complete.
\end{proof}

Then we give an application in degenerate Bourgain-Morrey-Sobolev spaces with matrix weights.
\begin{definition}
	Let $1 \le  p < t < r <\infty$ or $ 1 \le p \le t < r= \infty $.  Let $W: \rn \to  M_{n}  (\mathbb C)$ be a matrix weight and set the scalar weight $\omega : = \| W \|$.
	We define the degenerate Bourgain-Morrey-Sobolev spaces  $ \mathcal W ^ {1,p,t,r} (W)  $ by the set of all Lebesgue measurable functions on $\rn$
	such that
	\begin{equation*}
		\| f \| _{ \mathcal W ^ {1,p,t,r} (W) } := \| f \| _{ M^{t,r}_p ( \omega ) }   + \| \nabla f \| _{ M^{t,r}_p (W) } <\infty
	\end{equation*}
	where $\nabla f = ( \partial_{x_1} f , \ldots, \partial_{x_n}  f )^T$ is the gradient of $f$ and
	\begin{equation*}
		\| f \| _{ M^{t,r}_p ( \omega  ) } := \left\|     \left\{    \omega ( Q_{j,k} ) ^{ 1/t-1/p }   \left(  \int_{Q_{j,k}  }  |  f (y) |^p \omega(y) \d y     \right) ^{1/p}    \right\}    _{ j\in \mathbb Z , k\in \mathbb Z ^n }   \right\|   _{ \ell^r } ,
	\end{equation*}
	is  the norm of the scalar weighted Bourgain-Morrey space in (\ref{scalar wei BM}).
\end{definition}

\begin{lemma}[p. 262, \cite{W04}] \label{totally bounded}
	A set $\mathcal F$ in metric space $X$ is totally bounded if and only if  it is Cauchy-precompact, that is, every sequence has a Cauchy subsequence.
\end{lemma}

\begin{corollary}\label{cor Sobolev}
	Let $1 \le  p < t < r <\infty$ or $ 1 \le p \le t < r= \infty $. Let $W: \rn \to  M_{n}  (\mathbb C) $ be a matrix weight.
	A subset $\mathcal F \subset  \mathcal W ^ {1,p,t,r} (W)  $ is totally bounded if the following conditions are valid:

	{\rm  (i)} $\mathcal F$ uniformly vanishes at infinity, that is,
	\begin{equation*}
		\lim_{R\to \infty} \sup _{\vec f \in  \mathcal F} 	\| \vec f \chi_{B^c (0,R)}  \|_{   \mathcal W ^ {1,p,t,r} (W)    }  =0;
	\end{equation*}
	
	{\rm  (ii)} $\mathcal F$ is equicontinuous, that is,
	\begin{equation*}
		\lim _{a \to 0} \sup _{\vec f \in  \mathcal F} \sup_{y\in B(0,a)} 	\| \vec f  - \tau_y \vec  f  \|_{  \mathcal W ^ {1,p,t,r} (W)  } =0.
	\end{equation*}
\end{corollary}
\begin{proof}
	Note that $\mathcal F \subset  \mathcal W ^ {1,p,t,r} (W)  $ satisfies  conditions (i)-(ii) if and only if $\mathcal F \subset  M^{t,r}_p ( \omega ) $ satisfies conditions (i)-(ii) of Theorem \ref{K-R 1} ($d=1$) and  $\nabla \mathcal F : = \{ \nabla f : f\in \mathcal F \} \subset M^{t,r}_p (W) $ satisfies  (i)-(ii) of Theorem \ref{K-R 1} ($d=n$). Hence by Theorem \ref{K-R 1}, we obtain that both $ \mathcal F \subset  M^{t,r}_p ( \omega )  $ and $ \nabla \mathcal F \subset M^{t,r}_p (W)  $ are totally bounded. Then Corollary \ref{cor Sobolev} follows from Lemma \ref{totally bounded}.
\end{proof}

\section{Dyadic average operator conditions for precompact sets}

In this section, we give a criterion for precompactness by  dyadic average operator.
Let $\mathcal F$ be a subset of $M_p^{t,r} (W)$.	 $\mathcal F$ is equicontinuous  by means of the dyadic average operator if
\begin{equation}\label{aver con}
	\lim _{ a \to -\infty, a \in \mathbb Z } \sup _{\vec f \in  \mathcal F}  	\| \vec f  - E_{d,a}  \vec  f  \|_{ M_{p}^{t,r} (W)  } =0,
\end{equation}
where
\begin{equation} \label{Es}
	E_{d,a}  \vec f (x) := \sum_{Q \in  \D_{-a}}\frac{\chi_Q (x)}{|Q |} \int_{  Q  } \vec f (y) \d y.
\end{equation}
Recall that $x_Q$ is the lower left corner of $Q\in \D $.
Then for $x \in Q \in \mathcal D_{-a}$
\begin{equation*}
	E_{d,a}  \vec f (x) = \frac{\chi_{Q_{-a,  0}} (x -x_Q)}{|Q_{-a,  0} |} \int_{  Q_{-a,  0}  } \vec f (y -x_Q) \d y.
\end{equation*}
Then we will prove that (\ref{trans con}) is stronger than (\ref{aver con}). Indeed, suppose that condition (\ref{trans con}) holds. For any $\epsilon >0$, there exists a cube $R$ with center $0$ and side length $2^{a+1}>0$ such that
\begin{equation} \label{condin 3 trans}
	\sup _{\vec f \in  \mathcal F} \sup_{y\in R} 	\| \vec f  - \tau_y \vec  f  \|_{ M_{p}^{t,r} (W)  } <\epsilon.
\end{equation}
If $x,y \in Q \in \D_{-a}$, then $x-y \in R.$
By Jensen's inequality ($p\ge 1$), H\"older's inequality ($r/p \ge 1$), and (\ref{condin 3 trans}), we obtain
\begin{align*}
	&  \| \vec f - E_{d,a} \vec f \|_{  M_{p}^{t,r} (W)   } \\
	&  = \left(   \sum_{Q\in \D}W ( Q ) ^{r/t - r/p} \left( \int_Q    \left|  \frac{1}{| Q_{-a,  0} |}
	\int_{Q_{-a,  0}} W^{1/p}(x) \left ( \vec f (x) -  \chi_{Q_{-a,  0}} (x -x_Q) \vec f (y - x_Q ) \right )    \d y    \right| ^p \d x    \right)^{r/p} \right)^{1/r} \\
	\nonumber
	&\le  \left(   \sum_{Q\in \D}W ( Q ) ^{r/t - r/p}   \left( \frac{1}{| Q_{-a,  0}|} \int_{ Q_{-a,  0}}
	\int_Q    \left|   W^{1/p}(x) \left( \vec f (x) -\chi_{Q_{-a,  0}} (x -x_Q) \vec f (y - x_Q )  \right)     \right| ^p \d x  \d y      \right)^{r/p} \right)^{1/r} \\
	\nonumber
	&\le  \left(   \sum_{Q\in \D}W ( Q ) ^{r/t - r/p} \frac{1}{| Q_{-a,  0}|} \int_{ Q_{-a,  0}}
	\left(  \int_Q    \left|   W^{1/p}(x) \left( \vec f (x) - \chi_{Q_{-a,  0}} (x -x_Q) \vec f (y - x_Q )  \right)     \right| ^p \d x   \right)^{r/p} \d y      \right)^{1/r} \\
	\nonumber
	& \le \frac{1}{| Q_{-a,  0}|^{1/r}} | Q_{-a,  0}|^{1/r} \sup_{y\in R}  \left(   \sum_{Q\in \D}W ( Q ) ^{r/t - r/p}
	\left(  \int_Q    \left|   W^{1/p}(x) ( \vec f (x) - \vec f (x-y)  )     \right| ^p \d x   \right)^{r/p}      \right)^{1/r}  \le \epsilon,
\end{align*}
modified when $ r= \infty $.
Hence we prove that  (\ref{trans con}) is stronger than (\ref{aver con}).

Replacing the translation operator by the average operator, we have the following result.

\begin{theorem} \label{K-R 2}
	Let $1 \le  p < t < r <\infty$ or $ 1 \le p \le t < r= \infty $. Let $W :\rn \to  M_d (\mathbb C) $ be a matrix weight.
	A subset $\mathcal F \subset M_{p}^{t,r} (W) $ is totally bounded if the following conditions are valid:
	
	{\rm  (i)} $\mathcal F$ is bounded, that is,
	\begin{equation*}
		\sup _{\vec f \in  \mathcal F} 	\| \vec f \|_{  M_{p}^{t,r} (W)   } <\infty;
	\end{equation*}
	
	{\rm  (ii)} $\mathcal F$ uniformly vanishes at infinity, that is,
	\begin{equation*}
		\lim_{R\to \infty} \sup _{\vec f \in  \mathcal F} 	\| \vec f \chi_{B^c (0,R)}  \|_{  M_{p}^{t,r} (W)    }  =0;
	\end{equation*}
	
	{\rm  (iii)} $\mathcal F$ is equicontinuous by means of the dyadic average operator, that is, for any
	\begin{equation*}
		\lim _{ a \to -\infty, a \in \mathbb Z } \sup _{\vec f \in  \mathcal F}  	\| \vec f  - E_{d,a}  \vec  f  \|_{ M_{p}^{t,r} (W)  } =0.
	\end{equation*}
	where  $ E_{d,a}$ is same as (\ref{Es}).
\end{theorem}

\begin{proof}
	Assume that $\mathcal F \subset  M_{p}^{t,r} (W) $ satisfies (i)-(iii). Given $\epsilon >0$ small enough, to prove the total boundedness of $\mathcal F$, it suffices to find a finite  $\epsilon$-net of $\mathcal F$. Denote by $R_i := [-2^i , 2^i )^n$ for $i\in\mathbb Z$. Then from condition (ii), there exists a positive integer $m$ large enough such that
	\begin{equation} \label{epsilon 1 eva}
		\sup _{\vec f \in  \mathcal F} 	\| \vec f - \vec f  \chi_{R_m}  \|_{  M_{p}^{t,r} (W)   }  <\epsilon.
	\end{equation}
	Moreover, by condition (iii), there exists an integer $a<0 $ such that
	\begin{equation} \label{condi 2 eva}
		\sup _{\vec f \in  \mathcal F} \| \vec f  - E_{d,a}  \vec  f  \|_{ M_{p}^{t,r} (W)  } <\epsilon.
	\end{equation}
	There exists a sequence $\{ Q_j \} _{j=1}^N$ of disjoint cubes in $\D _{-a}$ such that $R_m = \bigcup_{i=1}^N  Q_j$, where $N=2^{  (m+1-a ) n}$.
	For any $\vec f \in \mathcal F$ and $x\in \rn$, let
	\begin{equation*}
		\Phi (\vec{ f}) (x):= \begin{cases}
			\vec f_{Q_j} :=\frac{1}{|Q_j|} \int_{Q_j} \vec f (y) \d y , &  x\in Q_j, j = 1,2,\ldots, N, \\
			\vec 0,  & \operatorname{otherwise}.
		\end{cases}
	\end{equation*}
	Note that for $x\in R_m$,
	\begin{equation*}
		\Phi (\vec{ f}) (x) = E_{d,a} \vec f (x).
	\end{equation*}
	Hence  via (\ref{epsilon 1 eva}) and (\ref{condi 2 eva}),  we have
	\begin{align}
		\nonumber
		\| \vec f - \Phi( \vec f) \|_{ M_{p}^{t,r} (W) } & \le 	\|  ( \vec f - \Phi( \vec f)  ) \chi_{R_m} \|_{ M_{p}^{t,r} (W) }  +	\|  ( \vec f - \Phi( \vec f)  ) \chi_{R_m^c} \|_{ M_{p}^{t,r} (W) } \\
		\nonumber
		& 
		=	\|  ( \vec f - E_{d,a} \vec f  ) \chi_{R_m} \|_{ M_{p}^{t,r} (W) }  +	\|   \vec f \chi_{R_m^c} \|_{ M_{p}^{t,r} (W) } \\
		& \le 		\|  ( \vec f - E_{d,a} \vec f  ) \|_{ M_{p}^{t,r} (W) } +\epsilon  \le 2\epsilon. \label{epsilon Phi eva}
	\end{align}
	From (\ref{epsilon Phi eva}), it suffices to show that $\Phi( \mathcal F)$ is totally bounded in $M_{p}^{t,r} (W)$. And we have proved that  $\Phi( \mathcal F)$ is totally bounded in $M_{p}^{t,r} (W)$ in Theorem \ref{K-R 1}. Thus the proof of Theorem \ref{K-R 2} is complete.
\end{proof}

\begin{lemma} [Lemma 4.5, \cite{CMR16}]  \label{matrix scalar}
	Let $1\le p <\infty$ and $W\in \mathcal A _p$. Then $ \|W\| $ and $ \|W^{-1} \| ^{-1} $ are scalar $A_p$ weights.
\end{lemma}

The following is a vector-valued extension of the Lebesgue differentiation theorem  on matrix weighted spaces $ M^{t,r}_p (W) $.
For any cube $Q$, $\vec f \in   L^1_{\operatorname{loc}} (\rn)  $, define
\begin{equation*}
	E_Q (\vec f) (x) = \frac{1}{|Q|} \int_Q  \vec f  (y) \d y.
\end{equation*}

\begin{theorem} \label{vec Les diff}
	Let $1 \le  p < t < r <\infty$ or $ 1 \le p \le t < r= \infty $. If $W\in \mathcal A _p $, then for any $\vec f \in M_{p}^{t,r} (W)$,
	\begin{equation*}
		\lim _{\ell(Q)  \to 0} \left| E_{ Q } \vec f (x) -\vec f (x)  \right| =0 \quad {\rm a.e.} \quad x \in \rn.
	\end{equation*}
\end{theorem}

\begin{proof}
	First, for any $\vec f = (f_1, f_2,\ldots, f_d)^T  \in M_{p}^{t,r} (W)$, it suffices to show that $f_i  \in L^1_{\operatorname{loc}} (\rn)$ for each $i=1,2,\ldots,d.$ Since $W\in \mathcal A _p$, then by  (\ref{A alpha})
	\begin{equation*}
		| \vec f |^p = |    W^{-1/p} W^{1/p} \vec f |^p \le \| W^{-1/p} \|^p   | W^{1/p} \vec f |^p = \| W^{-1} \|  | W^{1/p} \vec f |^p.
	\end{equation*}
	It follows that
	\begin{equation*}
		| \vec f |^p   \| W^{-1} \|^{-1} \le  | W^{1/p} \vec f |^p.
	\end{equation*}
	From Lemma \ref{matrix scalar}, we conclude that $  \| W^{-1} \|^{-1} $ is a scalar $A_p$ weight, hence $|\vec f | \in L^1_{\operatorname{loc}} (\rn)$. Indeed, since any compact set $K \subset \rn$ can be contained in a finite set of dyadic cubes, it suffices to show that for any $Q \in \D$, $| \vec f | \in L^1  (Q)$.
	By H\"older's inequality, we have
	\begin{align*}
		\nonumber
		\int_Q |\vec f (x)| \d x &= \int_Q  |\vec f(x)|    \| W^{-1} (x) \|^{-1/p}   \| W^{-1} (x) \|^{1/p}   \d x \\
		&\le  \left(  \int_Q  |\vec f (x)| ^p   \| W^{-1} (x)\|^{-1}  \d x \right)^{1/p}  \left(  \int_Q  \| W^{-1} (x) \|^{p'/p} \d x   \right)^{1/p'}  \\
		& \le C_{Q,W,p}    \left(  \int_Q | W^{1/p} (x) \vec f (x) |^p  \d x \right)^{1/p}
		\lesssim   \|\vec f \|_{M_{p}^{t,r} (W)   } <\infty.
	\end{align*}
	Hence, $f_i \in L^1_{\operatorname{loc}} (\rn)$ for each $i =1,2 ,\ldots, d$. By the classical Lebesgue differentiation theorem, for each $1 \le i \le d$, we have
	\begin{equation*}
		\lim _{\ell(Q)  \to 0} \left| E_Q  f_i (x) -f_i (x)  \right| =0 \quad {\rm a.e.} \quad x \in \rn.
	\end{equation*}
	Theorem \ref{vec Les diff} comes from the fact that for any $x\in \rn$,
	\begin{equation*}
		\left| E_Q \vec f (x) -\vec f (x)  \right|  \le d^{1/2} \max_i \left| E_Q  f_i (x) -f_i (x)  \right| .
	\end{equation*}
	Thus the proof is finished.
\end{proof}

Next we need the boundedness of dyadic average operator.
Given any collection $\mathcal Q$ of pairwise disjoint cubes $Q \subset
\rn$, define the averaging operator $A_{\mathcal Q}$ by
\begin{equation*}
	A_{\mathcal Q} \vec f (x) = \sum_{Q \in \mathcal Q} \frac{\chi_Q (x)}{|Q|} \int_Q \vec f (y) \d y .
\end{equation*}
\begin{lemma}[Proposition 4.7,\cite{CMR16}] \label{eva Lp}
	Let $1\le  p <\infty$.  Let $W:\rn \to  M_d (\mathbb C) $ be a matrix weight. Then $W \in \mathcal A _p$ if and only it satisfies
	\begin{equation*}
		\| 	A_{\mathcal Q} \vec f \|_{L^p (W) } \le C_{n,d,p,W} \| \vec f \|_{L^p (W)}.
	\end{equation*}
\end{lemma}

\begin{remark} \label{better estimate}
	In the sequel, let $p, p', W, \tilde d, \widetilde W , \tilde{\tilde d}$ have the same meaning as in Lemma \ref{impro AQ AR}.
	Let
	\begin{equation} \label{tilde beta}
		\tilde \beta := \begin{cases}
			n, & \operatorname{if} \; p =1, \\
			n +  \tilde{\tilde d} p / p',	& \operatorname{if} \;  1< p <\infty .
		\end{cases}
	\end{equation}
	Then we claim that for $ j+ a \in \mathbb N _0$,
	\begin{equation*}
		W ( (j + a)_{ \operatorname{pa} }Q  ) \lesssim \begin{cases}
			2^{(j+a)n} W(Q), &  \operatorname{if} p =1,\\
			2^{ (j+a) ( n +  \tilde{\tilde d} p / p') } W(Q), & \operatorname{if}  1< p <\infty .
		\end{cases}
	\end{equation*}
	Indeed, when $p =1$, by Lemma \ref{impro AQ AR}, we have
	\begin{align*}
		W ( (j + a)_{ \operatorname{pa} }Q  ) \approx |(j+a)_{\operatorname{pa}} Q| \| A_{ (j+a)_{\operatorname{pa}} Q  } \| \lesssim 2^{(j+a)n} |Q| \|A_Q\| \approx 2^{(j+a)n} W(Q).
	\end{align*}
	When $p \in (1,\infty)$, by Lemma \ref{impro AQ AR}, we obtain
	\begin{align*}
		W ( (j + a)_{ \operatorname{pa} }Q  ) & \approx | (j + a)_{ \operatorname{pa} }Q| \| A _{ (j + a)_{ \operatorname{pa} }Q  }  \|^p  \lesssim 2^{ (j+a) n} |Q| 2^{ (j+a)  \tilde{\tilde d} p / p'} \|A_Q\|^p  \approx 2^{ (j+a) (n +  \tilde{\tilde d} p / p' )  } W (Q).
	\end{align*}
	Thus we prove the claim. Now we show this estimate is better than $ 	W ( (j + a)_{ \operatorname{pa} }Q  ) \lesssim 2^{ (j+a) \beta} W (Q) $. Indeed, when $p = 1$, it is obvious since $\beta \ge n$. If $ p >1$,
	from \cite[Lemma 2.11 and Corollary 2.16]{BHYY23}, for $i \in \mathbb N _0$, we deduce that
	\begin{align*}
		\|A _{i _{ \operatorname{pa} } Q } A_Q ^{-1} \|^p & \approx \frac{1}{|i _{ \operatorname{pa} } Q | } \int_{  i _{ \operatorname{pa} } Q }  \left\| W ^{1/p} (x) A_Q ^{-1} \right\| ^p \d x 
		\\ &
		\approx \frac{1}{|i _{ \operatorname{pa} } Q | } \int_{  i _{ \operatorname{pa} } Q }  \left(  \frac{1}{|Q|} \int_Q \left\| W ^{1/p} (x) W ^{-1/p} (y) \right\|^{p'}  \d y  \right) ^{p/p'} \d x .
	\end{align*}
	From \cite[Lemma 2.2]{FR21}, we have  $\|A _{i _{ \operatorname{pa} } Q } A_Q ^{-1} \|^p  \lesssim 2^{i (\beta -n)}$. Hence
	\begin{equation*}
		\frac{1}{|i _{ \operatorname{pa} } Q | } \int_{  i _{ \operatorname{pa} } Q }  \left(  \frac{1}{|Q|} \int_Q \left\| W ^{1/p} (x) W ^{-1/p} (y) \right\|^{p'}  \d y  \right) ^{p/p'} \d x   \lesssim 2^{i (\beta -n)}.
	\end{equation*}
	This, together with \cite[Proposition 2.28(ii)]{BHYY23}, further implies that $W ^{ -1/(p-1)} \in \mathcal A _{p'}$ has
	the $\mathcal A _{p'}$-dimension $\tilde{\tilde d }: = (\beta -n) / (p-1)  $ and hence
	\begin{equation*}
		n + \tilde{\tilde d} p / p' = \beta .
	\end{equation*}
	Hence we  show this estimate is better than $ 	W ( (j + a)_{ \operatorname{pa} }Q  ) \lesssim 2^{ (j+a) \beta} W (Q) $.
\end{remark}

\begin{theorem} \label{tm ave BM}
	Let $W \in \mathcal A _p $ has the $\mathcal A _p$-dimension $\tilde d \in [0,n)$.
	Let $\tilde \beta$ be the same with (\ref{tilde beta}).
	
	{\rm (i)} Let $1 \le p <t< r < \infty $ .    If $ - nr/p +\tilde d r /p +n - \tilde \beta r (1/t-1/p) <0 $, then  for each $a \in \mathbb Z$, $ E_{d,a} $ is a bounded operator on matrix weighted Bourgain-Morrey spaces $M_p^{t,r} (W) $.
	
	{\rm (ii)} Let $1 \le p \le t< r = \infty $ .    If  $ \tilde d /p - n/p - \tilde \beta (1/t-1/p) \le 0$, then  for each $a \in \mathbb Z$, $ E_{d,a} $ is a bounded operator on matrix weighted Bourgain-Morrey spaces $M_p^{t,r} (W) $.
\end{theorem}
\begin{proof}
	(i) For a fixed $a\in \mathbb Z$, we have
	\begin{align*}
		\| E_{d,a} \vec f \|_ {M_p^{t,r} (W) } & \le \left(   \sum_{j \le -a} \sum_{Q\in \mathcal D_j} W(Q)^{r/t-r/p}  \left( \int_Q     | W^{1/p}(x) E_{d,a}  \vec f (x) | ^p \d x    \right)^{r/p} \right)^{1/r}  \\
		& + \left(   \sum_{j > -a} \sum_{Q\in \mathcal D_j} W(Q)^{r/t-r/p}  \left( \int_Q     | W^{1/p}(x) E_{d,a}  \vec f (x) | ^p \d x    \right)^{r/p} \right)^{1/r}   =: S_1 +S_2.
	\end{align*}
	By Lemma \ref{eva Lp}, we obtain
	\begin{align*}
		S_1 &\le  C_{n,d,p,W} \left(   \sum_{j \le -a} \sum_{Q\in \mathcal D_j} W(Q)^{r/t-r/p}  \left( \int_Q     | W^{1/p}(x)   \vec f (x) | ^p \d x    \right)^{r/p} \right)^{1/r}
		\lesssim \| \vec f \|_ {M_p^{t,r} (W) } .
	\end{align*}
	When $j > -a$, we denote by $ (j + a)_{ \operatorname{pa} }Q \in \D _{-a}$ the $(j + a)$-th dyadic parent of $Q \in \D _j $.
	Let $ \{A_Q\} _{Q \in \D} $ be the reducing operator of order $p$ for $W$.
	By Lemmas \ref{impro AQ AR},  we have
	\begin{equation} \label{AQ LE}
		\| A_Q A_{ (j + a)_{ \operatorname{pa} }Q  }^{-1} \| \le c 2^{ (j+a)\tilde d /p } .
	\end{equation}
	In $S_2 $, using (\ref{reducing A_Q}) and  (\ref{AQ LE})  , we have
	\begin{align*}
		S_2^r &= \sum_{j > -a} \sum_{Q\in \mathcal D_j} W(Q)^{r/t-r/p}  \left( \int_Q     | W^{1/p}(x) E_{d,a}  \vec f (x) | ^p \d x    \right)^{r/p} \\
		& = \sum_{j > -a} \sum_{Q\in \mathcal D_j} W(Q)^{r/t-r/p}
		\left( \frac{|Q|}{|Q|} \int_Q    \left| W^{1/p}(x) \frac{1}{|(j + a)_{ \operatorname{pa} }Q | } \int_{ (j + a)_{ \operatorname{pa} }Q } \vec f (y)  \d y  \right| ^p \d x    \right)^{r/p} \\
		& \approx  \sum_{j > -a} \sum_{Q\in \mathcal D_j}  W(Q)^{r/t-r/p} |Q|^{r/p}  \left(  \left|A_Q  \frac{1}{|(j + a)_{ \operatorname{pa} }Q | } \int_{ (j + a)_{ \operatorname{pa} }Q } \vec f (y)  \d y  \right|   \right)^{r} \\
		& \lesssim \sum_{j > -a} \sum_{Q\in \mathcal D_j} W(Q)^{r/t-r/p} |Q|^{r/p}
		\left( 2^{ (j+a)\tilde d /p } \left|A_{ (j + a)_{ \operatorname{pa} }Q  }  \frac{1}{|(j + a)_{ \operatorname{pa} }Q | } \int_{ (j + a)_{ \operatorname{pa} }Q } \vec f (y)  \d y  \right|   \right)^{r} \\
		& \approx \sum_{j > -a} \sum_{Q\in \mathcal D_j}W(Q)^{r/t-r/p} |Q|^{r/p}  2^{ (j+a)\tilde d r/p }
		\left(  \frac{1}{| (j + a)_{ \operatorname{pa} }Q  | } \int_{ (j + a)_{ \operatorname{pa} }Q  }     | W^{1/p}(x) \frac{1}{|(j + a)_{ \operatorname{pa} }Q | } \int_{ (j + a)_{ \operatorname{pa} }Q } \vec f (y)  \d y  | ^p \d x    \right)^{r/p}  \\
		& = \sum_{j > -a} \sum_{Q\in \mathcal D_j} W(Q)^{r/t-r/p} 2^{-jnr/p}  2^{ (j+a)\tilde d r/p }
		\left(  \frac{1}{| (j + a)_{ \operatorname{pa} }Q  | } \int_{ (j + a)_{ \operatorname{pa} }Q  }     | W^{1/p}(x) E_{d,a} \vec f (x)  | ^p \d x    \right)^{r/p}  \\
		& = \sum_{j > -a} \sum_{Q\in \mathcal D_j} W(Q)^{r/t-r/p} 2^{-jnr/p}  2^{ (j+a)\tilde d r/p } 2^{-anr/p}
		\left(  \int_{ (j + a)_{ \operatorname{pa} }Q  }     | W^{1/p}(x) E_{d,a} \vec f (x)  | ^p \d x    \right)^{r/p} .
	\end{align*}
	Remark that given $S \in \D _{-a} $, there are $ 2^{ (j+a )n} $ cubes $R $ such that $  (j+a)_{\operatorname{pa}} R = S$.
	From  Remark \ref{better estimate} and  $1/t - 1/p <0$, we have
	\begin{equation*}
		W ( Q  ) ^{1/t - 1/p} \lesssim \begin{cases}
			2^{-(j+a)n ( 1/t - 1)} W( (j + a)_{ \operatorname{pa} }Q)^{1/t - 1}, &  \operatorname{if} p =1,\\
			2^{ - (j+a)(n +  \tilde{\tilde d} p / p' )  (1/t - 1/p) } W( (j + a)_{ \operatorname{pa} }Q)^{1/t - 1/p}, & \operatorname{if}  1< p <\infty .
		\end{cases}
	\end{equation*}
	That is  $W ( Q  ) ^{1/t - 1/p} \lesssim  	2^{ - (j+a) \tilde \beta (1/t - 1/p) } W( (j + a)_{ \operatorname{pa} }Q)^{1/t - 1/p} $.
	Hence, by Lemma \ref{eva Lp} and  $ - nr/p +\tilde d r /p +n - \tilde \beta r (1/t-1/p) <0 $,  we have
	
	\begin{align*}
		& \sum_{j > -a} \sum_{Q\in \mathcal D_j} W(Q)^{r/t-r/p} 2^{-jnr/p}  2^{ (j+a)\tilde d r/p } 2^{-anr/p}
		\left(  \int_{ (j + a)_{ \operatorname{pa} }Q  }     | W^{1/p}(x) E_{d,a} \vec f (x)  | ^p \d x    \right)^{r/p} \\
		& =	 \sum_{j > -a}2^{-jnr/p}  2^{ (j+a)\tilde d r/p } 2^{-anr/p}
		\sum_{Q\in \mathcal D_j} W(Q)^{r/t-r/p}    \left(  \int_{ (j + a)_{ \operatorname{pa} }Q  }     | W^{1/p}(x) E_{d,a} \vec f (x)  | ^p \d x    \right)^{r/p} \\
		& \le \sum_{j > -a}2^{-jnr/p}  2^{ (j+a)\tilde d r/p } 2^{-anr/p}
		2^{(j+a)n  } 2^{-(j+a)  \tilde \beta (r/t-r/p)}
		\sum_{S\in \mathcal D_{-a}} W(S)^{r/t-r/p}    \left(  \int_{ S }     | W^{1/p}(x) E_{d,a} \vec f (x)  | ^p \d x    \right)^{r/p} \\
		& \lesssim  \| \vec f \|_{M_p^{t,r} (W)   } ^r.
	\end{align*}

	(ii) Fix $Q_{j,k} \in \D$. If $ j \le -a$, then by Lemma \ref{eva Lp}, we have
	\begin{align*}
		&	W ( Q_{j,k} ) ^{1/t-1/p} \left(  \int_{ Q_{j,k}   } |W^{1/p} (x)  E_{d,a} \vec f (x)  | ^p \d x    \right)^{1/p} \lesssim 	W ( Q_{j,k} ) ^{1/t-1/p} \left(  \int_{ Q_{j,k}   } |W^{1/p} (x) \vec f (x)  | ^p \d x    \right)^{1/p}.
	\end{align*}
	If $j> -a $, we denote by $ (j + a)_{ \operatorname{pa} } Q_{j,k}  \in \D _{-a}$ the $(j + a)$-th dyadic parent of $Q_{j,k} \in \D _j $.
	\begin{align*}
		&	W ( Q_{j,k} ) ^{1/t-1/p} \left(  \int_{ Q_{j,k}   } |W^{1/p} (x)  E_{d,a} \vec f (x)  | ^p \d x    \right)^{1/p}\\
		& = 	W ( Q_{j,k} ) ^{1/t-1/p}  \left( \frac{|Q_{j,k}|}{|Q_{j,k}|} \int_{ Q_{j,k} }    \left| W^{1/p}(x) \frac{1}{|(j + a)_{ \operatorname{pa} } Q_{j,k} | } \int_{ (j + a)_{ \operatorname{pa} } Q_{j,k} } \vec f (y)  \d y  \right| ^p \d x \right)^{1/p} \\
		& \approx 	W ( Q_{j,k} ) ^{1/t-1/p} |Q_{j,k}|^{1/p}    \left| A_{Q_{j,k}  } \frac{1}{|(j + a)_{ \operatorname{pa} } Q_{j,k} | } \int_{ (j + a)_{ \operatorname{pa} } Q_{j,k} } \vec f (y)  \d y \right|  \\
		& \lesssim	W ( Q_{j,k} ) ^{1/t-1/p} |Q_{j,k}|^{1/p}  2^{(j+a) \tilde d /p }  \left| A_{(j + a)_{ \operatorname{pa} } Q_{j,k}   } \frac{1}{|(j + a)_{ \operatorname{pa} } Q_{j,k} | } \int_{ (j + a)_{ \operatorname{pa} } Q_{j,k} } \vec f (y)  \d y \right|  \\
		& \approx	W ( Q_{j,k} ) ^{1/t-1/p} 2^{-jn/p}  2^{(j+a) \tilde d /p }  	\left( \frac{1}{|(j + a)_{ \operatorname{pa} }Q _{j,k}|}
		\int_{ (j + a)_{ \operatorname{pa} } Q_{j,k}  }    \left| W^{1/p}(x) \frac{1}{|(j + a)_{ \operatorname{pa} } Q_{j,k} | } \int_{ (j + a)_{ \operatorname{pa} } Q_{j,k} } \vec f (y)  \d y  \right| ^p \d x \right)^{1/p} .
	\end{align*}
	Then by Lemma \ref{eva Lp} and $ \tilde d /p - n/p - \tilde \beta (1/t-1/p) \le 0$,   we obtain
	\begin{align*}
		&	W ( Q_{j,k} ) ^{1/t-1/p} \left(  \int_{ Q_{j,k}   } |W^{1/p} (x)  E_{d,a} \vec f (x)  | ^p \d x    \right)^{1/p} \\
		& \lesssim W ( Q_{j,k} ) ^{1/t-1/p} 2^{-jn/p}  2^{(j+a) \tilde d /p }
		\left( \frac{1}{|(j + a)_{ \operatorname{pa} }Q _{j,k}|} \int_{ (j + a)_{ \operatorname{pa} } Q_{j,k}  }    \left| W^{1/p}(x) \vec f (x) \right|^p \d  x \right)^{1/p} \\
		& = W ( Q_{j,k} ) ^{1/t-1/p} 2^{-jn/p}  2^{(j+a) \tilde d /p } 2^{-an/p} \left(\int_{ (j + a)_{ \operatorname{pa} } Q_{j,k}  }    \left| W^{1/p}(x) \vec f (x) \right|^p \d  x \right)^{1/p} \\
		& \lesssim  2^{(j+a) (\tilde d /p -n/p) } 2^{- (j+a) \tilde \beta (1/t-1/p)} W ( (j + a)_{ \operatorname{pa} }Q _{j,k} ) ^{1/t-1/p}
		\left(\int_{ (j + a)_{ \operatorname{pa} } Q_{j,k}  }
		\left| W^{1/p}(x) \vec f (x) \right|^p \d  x \right)^{1/p} \\
		& \lesssim  W ( (j + a)_{ \operatorname{pa} }Q _{j,k} ) ^{1/t-1/p} \left(\int_{ (j + a)_{ \operatorname{pa} } Q_{j,k}  }    \left| W^{1/p}(x) \vec f (x) \right|^p \d  x \right)^{1/p} .
	\end{align*}
	Finally, take the supremum over the cubes $Q_{j,k} \in D$, we obtain
	\begin{equation*}
		\| E_{d,a} \vec f \|_{ M_p^{t,\infty} (W) }  \lesssim \|  \vec f \|_{ M_p^{t,\infty} (W) } .
	\end{equation*}
	Hence we finish the proof.
\end{proof}

Finally, we have the following Kolmogorov-Riesz compactness theorem for matrix weighted Bourgain-Morrey spaces.
\begin{theorem}\label{char Bourgain}
	Let $1 \le  p  <t <r<\infty$ . Let $W \in \mathcal A _p $ with the $\mathcal A _p$-dimension $\tilde d \in [0,n)$.
	Let $\tilde \beta$ be the same with (\ref{tilde beta}).
	Let   $  - nr/p +\tilde d r /p +n - \tilde \beta r (1/t-1/p) <0  $. A subset $\mathcal F$ of $M_{p}^{t,r} (W) $ is totally bounded if and only if the following conditions hold:
	
	{\rm  (i)} $\mathcal F$ is bounded, that is,
	\begin{equation*}
		\sup _{\vec f \in  \mathcal F} 	\| \vec f \|_{  M_{p}^{t,r} (W)   } <\infty;
	\end{equation*}
	
	{\rm  (ii)} $\mathcal F$ uniformly vanishes at infinity, that is,
	\begin{equation*}
		\lim_{R\to \infty} \sup _{\vec f \in  \mathcal F} 	\| \vec f \chi_{B^c (0,R)}  \|_{  M_{p}^{t,r} (W)    }  =0;
	\end{equation*}
	
	{\rm  (iii)} $\mathcal F$ is equicontinuous by means of the dyadic average operator, that is, for any
	\begin{equation*}
		\lim _{ a \to -\infty, a \in \mathbb Z } \sup _{\vec f \in  \mathcal F}  	\| \vec f  - E_{d,a}  \vec  f  \|_{ M_{p}^{t,r} (W)  } =0.
	\end{equation*}
	where  $ E_{d,a}$ is same as (\ref{Es}).

\end{theorem}

\begin{proof}
	The sufficiency is due to Theorem \ref{K-R 2}. Now we prove the necessity.
	Assume that $\mathcal F$ of $M_{p}^{t,r} (W) $ is totally bounded.
	For any given $\epsilon >0$, there exists $ \{  \vec f _k \} _{k=1}^{N_0} \subset \mathcal F$ such that  $ \{  \vec f _k \} _{k=1}^{N_0}$ is an $\epsilon$-net of $\mathcal F$, that is, for any $\vec f \in \mathcal F$, there exists $\vec f _k, k \in \{1,2,\ldots, N_0\} $ such that $  \| \vec f - \vec f_k  \|_{ M_{p}^{t,r} (W) } < \epsilon $.
	
	Clearly, (i) is true.
	
	For any $\epsilon>0$ and $k \in\{1,2,\ldots, N_0\} $, since $r<\infty$, there exists $R_k>0$ such that
	\begin{equation*}
		\| \vec f _k \chi_{ B^c (0,R_k) }  \|_{ M_{p}^{t,r} (W)  } <\epsilon.
	\end{equation*}
	Taking $R =\max_{  k \in \{1,\ldots, N_0 \}  } \{ R_k \}$, we have $	\| \vec f _k \chi_{ B^c (0,R) }  \|_{ M_{p}^{t,r} (W)  } <\epsilon.$
	Then for any $\vec f\in \mathcal F$,
	\begin{align*}
		\| \vec f \chi_{B^c (0,R)}  \|_{  M_{p}^{t,r} (W)    }  \le 	\| ( \vec f -\vec f_k )\chi_{B^c (0,R)}  \|_{  M_{p}^{t,r} (W)    }  +	\| \vec f_k \chi_{B^c (0,R)}  \|_{  M_{p}^{t,r} (W)    }   < 2\epsilon.
	\end{align*}
	Thus we prove that
	\begin{equation*}
		\lim_{R\to \infty} \sup _{\vec f \in  \mathcal F} 	\| \vec f \chi_{B^c (0,R)}  \|_{  M_{p}^{t,r} (W)    }  =0.
	\end{equation*}
	Therefore, we prove that $\mathcal F$ satisfies condition (ii).

	As for (iii), for each $ 1\le k \le N_0$,  by Theorem \ref{tm ave BM}, we have for each $a \in \mathbb Z$,
	\begin{equation*}
		\|  \vec f_k  - E_{d,a}  \vec  f_k \|_ {  M_{p}^{t,r} (W) }   \le (1+c) 	\|  \vec f_k \|_ {  M_{p}^{t,r} (W) }  .
	\end{equation*}
	Thus using the   dominated convergence theorem   to  obtain that there exists $ a_0\in \mathbb Z $ such that for any $a \le  a_0$,
	\begin{align*}
		\max_{1\le  k \le N_0} 	\|  \vec f _k - E_{d,a}  \vec  f _k \|_ {  M_{p}^{t,r} (W) } <\epsilon.
	\end{align*}
	Now if $a \le a_0$, then
	\begin{align*}
		\|  \vec f  - E_{d,a}  \vec  f \|_ {  M_{p}^{t,r} (W) } & \le 	\| E_{d,a} \vec f  - E_{d,a}  \vec  f_k \|_ {  M_{p}^{t,r} (W) } + 	\| E_{d,a} \vec f_k  -   \vec  f_k \|_ {  M_{p}^{t,r} (W) }  + 	\|  \vec f_k  -  \vec  f \|_ {  M_{p}^{t,r} (W) }  \\
		& \le c 	\|  \vec f  -  \vec  f_k \|_ {  M_{p}^{t,r} (W) } + \epsilon + \epsilon
		\lesssim \epsilon.
	\end{align*}
	Then the proof is complete.
\end{proof}

\begin{remark}
	Theorem \ref{char Bourgain} is not true for $ 1< p<t<r = \infty   $.
	Indeed, let $d=1, W\equiv 1$, $0<p<t<r = \infty  $. Let $f (x) = |x|^{-n/t}$. Then
	\begin{align*}
		\| f \| _{ M ^{t,\infty}_{p}   } &
		= \sup_{Q \in \D} |Q|^{1/t-1/p} \left(  \int_Q |y|^{ -np /t  } \d y \right)^{1/p}
		\approx  1.
	\end{align*}
	Let $f_j  = \chi_{B(0,2^j)} f $ for $j \in \mathbb N$. Note that $f$ and $f_j$ are  radial and symmetric functions. Let  $x = ( x_1, 0,\ldots,0 )$ where $x_1  = 2^j  + 2\times 2^j$. Let  $s = x_1 /10  $ . Then $ B(x, s)  \subset B(0,2^j)^c $.
	Then
	\begin{align*}
		\| f - f_j \| _{   M ^{t,\infty}_{p}  }
		& \ge c s^{ -n (1/t-1/p) }  \left(  \int_{B(x,s)}   |y| ^{  -np /t }  \d y   \right) ^{1/p}  \ge c  s^{ -n (1/t-1/p) } \left(  \int_{|x| -|s|} ^{ |x| +|s| }  \eta^{-np/t } \eta^{n-1} \d \eta         \right)^{1/p} \\
		& =c \left( \frac{|x|}{10} \right)^{ -n (1/t-1/p) }  \left( \left(1+\frac{1}{10}\right)^{n-np/t} - \left(1-\frac{1}{10}\right)^{n-np/t}    \right)^{1/p} |x|^{ n (1/t-1/p ) } \\
		&= c >0.
	\end{align*}
	This shows that the totally bounded set $\mathcal F = \{ f \}$ of $M ^{t,\infty}_{p}$ is not uniformly vanishes at infinity.
\end{remark}

\end{document}